\documentclass[11pt]{article}
\usepackage[T2A]{fontenc}
\usepackage[english,russian]{babel}
\usepackage[tbtags]{amsmath}
\usepackage{amsfonts,amssymb,mathrsfs,amscd}

\usepackage[utf8]{inputenc} 

\usepackage{geometry} 
\usepackage{graphicx} 

\usepackage{booktabs} 
\usepackage{array} 
\usepackage{paralist} 
\usepackage{verbatim} 
\usepackage{subfig} 

\newtheorem{theorem}{Теорема}[section]
\newtheorem{lemma}{Лемма}[section]
\newtheorem{corollary}{Следствие}[section]

\newtheorem{definition}{Определение}[section]

\newtheorem{remark}{Замечание}[section]

\begin{document}

\title{Ergodicity conditions for general Markov chains in terms of invariant finitely additive measures}
\author{Alexander ~I.~Zhdanok \footnote{Institute for Information Transmission Problems of the Russian Academy of Science, Moscow, Russia;} \footnote{Tuvinian Institute for Exploration of Natural Resources of the Siberian Branch RAS, Kyzyl, Republic of Tuva, Russia; e-mail: zhdanok@inbox.ru;} \footnote{This work was supported by the Russian Foundation for Basic Research (grant № 20-01-00575-a)}}

\date{}

\maketitle

{\bf Abstract:} We consider general Markov chains with discrete time in an arbitrary measurable (phase) space and homogeneous in time. Markov chains are defined by the classical transition function which within the framework of the operator treatment generates a conjugate pair of linear Markov operators in the Banach space of measurable bounded functions and in the Banach space of bounded finite additive measures. It is proved that the well-known Doeblin condition $ (D) $ of ergodicity (quasi\-compactness) of the Markov chain is equivalent to the condition $ (*) $: all finitely additive invariant measures of the Markov operator are countably additive  i.e. there are no invariant purely finitely additive measures. Under some assumptions, it is proved that the conditions $ (D) $ and $ (*) $ are also equivalent to the condition $ (**) $: the set of invariant finitely additive measures of a Markov operator is finite-dimensional. Ergodic theorems are given. \\



\newpage
\begin{center}
{\bf Условия эргодичности общих цепей Маркова \\ в терминах инвариантных конечно-аддитивных мер }\\
\end{center}
\maketitle
\begin{center}
{А.~И.~Жданок
\footnote{Институт проблем передачи информации им.~А.\,А.~Харкевича Российской академии наук, г. Москва, Россия;} 
\footnote{Тувинский институт комплексного освоения природных ресурсов Сибирского отделения Российской академии наук, г. Кызыл, Россия; e-mail: zhdanok@inbox.ru;}
\footnote{Работа выполнена при поддержке Российского фонда фундаментальных исследований (грант No~20-01-00575-а).}}
\end{center}


{\bf Аннотация:} Рассматриваются общие цепи Маркова с дискретным временем в произвольном измеримом (фазовом) пространстве и однородные по времени. Цепи Маркова задаются классической переходной функцией, которая в рамках операторного подхода порождает сопряженную пару линейных марковских операторов в банаховом пространстве измеримых ограниченных функций и в банаховом пространстве ограниченных конечно-ад\-ди\-тив\-ных мер. Доказывается, что известное условие Деблина $(D)$ эргодичности (квазикомпактности) цепи Маркова эквивалентно условию $(*)$: все конечно-аддитивные инвариантные меры марковского оператора являются счётно-ад\-ди\-тив\-ными, т. е. отсутствуют инвариантные чисто конечно-аддитивные меры. При некоторых допущениях доказывается, что условия $(D)$ и $(*)$ эквивалентны также условию $(**)$: множество инвариантных конечно-ад\-ди\-тив\-ных мер марковского оператора конечномерно. Даются эргодические теоремы. 

{\bf Ключевые слова:} общие цепи Маркова, марковские операторы, конеч\-но-ад\-ди\-тив\-ные меры, инвариантные меры, условия ква\-зиком\-пакт\-нос\-ти, эргодические теоремы.

\tableofcontents

\newpage
{\bf Предисловие}

80 лет тому назад, в 1940 году, в журнале ``Математический сборник''
была опубликована (на англ. яз.) первая часть фундаментальной статьи А. Д. Александрова ``Additive set-functions in abstract spaces'' ~\cite{Ale1}, а затем в 1941 и в 1943 годах были там же опубликованы вторая и третья её части ~\cite{Ale2}, ~\cite{Ale3}. Эта очень большая работа (суммарный объем -- 170 журнальных страниц) была первой в мировой литературе, целиком посвящённой построению основ общей Теории конечно-аддитивных мер (``зарядов'') в топологических пространствах. Она опиралась на предыдущие статьи различных авторов, проводивших исследования по частным вопросам этой новой теории, но содержала и большое количество принципиально новых (на тот период) конструкций и теорем, принадлежащих А. Д. Александрову.
Эта работа А. Д. Александрова востребована и сегодня.

В своих работах ~\cite{Zhd01}, ~\cite{Zhd02} и в других мы также использовали некоторые результаты из \cite{Ale1}, ~\cite{Ale2},  ~\cite{Ale3}, и, прежде всего, ``самую красивую'' (по мнению многих) теорему А. Д. Александрова о том, что любая регулярная конечно-аддитивная мера на компакте является счётно-аддитивной (\cite[гл.~III, \S 9, теорема 5, стр. 590]{Ale2}). Эта и несколько других  теорем А. Д. Александрова приводятся и в монографии  ~\cite[главы III, IV]{DS1}, наиболее полно отражающей развитие теории конечно-аддитивных мер на конец 1950-ых годов со своими дополнениями.

А. Д. Александров, по-видимому, был первый, кто использовал компактификацию исходного пространства для изучения заданных на нем конечно-адди\-тив\-ных мер. За основу им была взята конструкция Уолмэновского расширения (1937 год), существенно затем развитая в ~\cite{Ale1}, ~\cite{Ale2}, ~\cite{Ale3}. Эта идея оказалась весьма плодотворной -- во многих дальнейших работах разных авторов по теории конечно\--ад\-ди\-тив\-ных мер строятся свои модели подобных компактификаций. Мы также в статьях ~\cite{Zhd01}, ~\cite{Zhd02}, а затем и в ~\cite{Zhd04} строим свою конструкцию, которую назвали ``Гамма-компак\-ти\-фи\-ка\-цией измеримых пространств'', и применяем её для исследования общих цепей Маркова в ~\cite{Zhd02}.

Ниже мы приведем свои комментарии по одному важному результату А. Д. Александрова, используемому нами в основном тексте настоящей статьи в $\S 1$ (и далее).\medskip

В настоящей статье мы рассматриваем общие цепи Маркова (ЦМ), понимаемые как случайный марковский процесс с дискретным временем в произвольном фазовом пространстве и однородный по времени.  Цепи Маркова задаются переходной вероятностью (функцией), счётно-аддитивной по второму аргументу, т. е. мы рассматриваем только ``классические'' цепи Маркова. Переходные вероятности порождают интегральные марковские операторы, действующие в пространствах измеримых функций, счётно-аддитивных и конечно-аддитивных мер. В рамках такого операторного подхода в статье в $\S 4$ и $\S 5$ получены следующие основные результаты.

В Теореме 4.4 доказывается, что известное условие Деблина $(D)$ эргодичности произвольной ЦМ эквивалентно введённому нами ранее условию $(*)$: все инвариантные конечно-аддитивные меры ЦМ являются счётно-аддитивными, что означает, что у ЦМ нет инвариантных чисто конечно\--ад\-ди\-тив\-ных мер (т.е. выполняется эквивалентное наше условие $(\tilde{*})$). Отмечается также, что из условия $(*)$ следует известное условие квазикомпактности ЦМ $(K)$.

Рассматривается еще одно наше условие $(**)$: множество инвариантных для ЦМ конечно-аддитивных мер имеет конечную размерность $n$. Ранее нами было доказано (здесь это Теорема 5.1), что выполнение условия $(*)$ влечет за собой выполнение условия $(**)$. В Теореме 5.3 доказывается обращение этого утверждения в одномерном случае $n=1$. Это означает, что если ЦМ имеет только одну инвариантную конечно-аддитивную меру, то она счётно-аддитивна и выполнены наше условие $(*)$, условие Деблина $(D)$ и условие квазикомпактности $(K)$ с соответствующими эргодическими следствиями.

В Теореме 5.5 доказывается обращение упомянутой Теоремы 5.1 уже для произвольной конечной размерности $n\ge2$, но при некоторых дополнительных условиях $(\alpha)$ и $(\beta)$. В ней утверждается, что из условия $(**)$ следует условие $(*)$, а также эквивалентное ему условие Деблина $(D)$, откуда следует условие квазикомпактности $(K)$. 

Как следствие из полученных результатов мы приводим равномерные эргодические Теорему 4.5, Теорему 5.3 (для $n=1$) и Теорему 5.6 (для $n\ge2$). 

Подчеркнем, что основные полученные результаты верны для любого измеримого фазового пространства ЦМ $(X, \Sigma)$, в том числе для счётных ЦМ, для дискретных $\Sigma=2^{X}$, для произвольных топологических пространств $X$ с борелевской сигма-алгеброй $\Sigma=\mathfrak B$, для $X=[0, 1]$, $X=R$ и т.д.

В основной теоретической части статьи ($\S 1$ -- $\S 5$) на переходную вероятность ЦМ никаких дополнительных требований, кроме перечисленных выше условий, не накладывается.

В $\S 1$ -- $\S 5$ статьи случайные величины (элементы), порождаемые цепью Маркова, не используются.

Представляемую здесь работу можно охарактеризовать и как применение общей теории конечно-аддитивных мер в операторной теории общих цепей Маркова.
\section{Определения, термины, обозначения и некоторые известные сведения}\label{s1}

В настоящей статье мы будем использовать некоторые работы и результаты по теории меры и функциональному анализу, начиная с 1930-ых годов. За прошедшее время символика и терминология в этих направлениях существенно изменилась, есть различия и в современной литературе. Поэтому приведём некоторые используемые нами базовые определения и понятия  и их символику, которыми сегодня чаще всего пользуются и другие авторы. Здесь и далее мы ориентируемся на определения и символику в ~\cite{DS1} и ~\cite{YoHew1}.

Пусть $X$ – произвольное бесконечное множество.

\begin{definition}\label{de1.1}
 Семейство подмножеств $\Sigma$ в множестве $X$ называется {\it алгеброй} его подмножеств, если $\O , X \in \Sigma$ и, из того, что $A, B \in \Sigma$ следует, что $A \cup B, A \cap B, A \setminus B \in \Sigma$.
 \end{definition}
 
\begin{definition}\label{de1.2}
Семейство подмножеств $\Sigma$ в множестве $X$ называется {\it сигма-алгеброй} ($\sigma$-алгеброй), если оно является алгеброй, и выполнено дополнительно условие: если $A_{1}, A_{2}, A_{3}, ... \in \Sigma$, то $\bigcup_{n=1}^{\infty} A_{n} \in \Sigma$ и $\bigcap_{n=1}^{\infty}A_{n} \in \Sigma$.
\end{definition}

\begin{remark}\label{r1}
Указанные в двух определениях условия сделаны для удобства избыточными. 
\end{remark}

\begin{definition}\label{de1.3}
 Пара $(X,\Sigma)$ называется {\it измеримым пространством}, а каждое множество $A \in \Sigma$ - {\it измеримым}.
 \end{definition}

Сегодня эти термины обычно вводятся до того, как на $(X,\Sigma)$ будет задана какая-либо мера, и даже когда она не будет задана.

В статье будут рассматриваться измеримые пространства $(X,\Sigma)$, где $\Sigma$ является сигма-алгеброй.

Всегда (по умолчанию) мы будем предполагать, что сигма-алгебра $\Sigma$ содержит все одноточечные множества  $\{x\} \in \Sigma, x \in X$. Это стандартное требование, даже когда оно явно не оговаривается.

Всюду далее $R=R^{1}$ - множество действительных чисел (числовая прямая).

Обозначим $B(X,\Sigma)$ - банахово пространство ограниченных $\Sigma$ -измеримых функций $f\colon X\to R$ с $sup$-нормой $\| f \|=sup |f(x)|$.

Любую функцию  $\varphi \colon \Sigma \to R$ будем называть функцией множеств на $\Sigma$.
Функция множеств $\varphi (E), E \in \Sigma$, может принимать как положительные, так и отрицательные значения. В настоящей работе мы используем только ограниченные функции множеств, т. е. такие, что $sup |\varphi(E)|<\infty$, где супремум берется по всем множествам $E \in \Sigma$.

 \begin{definition}\label{de1.4}
 Пусть $(X,\Sigma)$ произвольное измеримое пространство. Ограниченная функция множеств  $\mu \colon \Sigma \to R$ называется {\it конечно-аддитивной мерой}, если $\mu(\O)=0$, и для любых множеств $E_{1}, E_{2} \in \Sigma$ таких, что $E_{1} \cap E_{2}=\O$, выполняется $\mu (E_{1} \cup E_{2}) = \mu (E_{1}) + \mu (E_{2})$.
 \end{definition}
 
 Конечно-аддитивные меры в литературе также называют {\it аддитивными мерами, средними, зарядами} и др.

\begin{definition}\label{de1.5}
 Пусть $(X,\Sigma)$ произвольное измеримое пространство. Ограниченная функция множеств  $\mu \colon \Sigma \to R$ называется {\it счётно-аддитивной мерой}, если она является конечно-аддитивной мерой и выполнено условие: если $E_{1}, E_{2}, ... \in \Sigma$, $E_{i} \cap E_{j}=\O$ при $i \ne j$, то
 $$
 \mu (\bigcup_{n=1}^{\infty} E_{n}) = \sum_{n=1}^{\infty} \mu (E_{n}). 
 $$
 \end{definition}
  
 Счётно-аддитивные меры в литературе также называют {\it вполне аддитивными мерами}.

\begin{definition}\label{de1.6}
(см. ~\cite{YoHew1}). Конечно-аддитивная неотрицательная мера $\mu \colon \Sigma \to R$ называется {\it чис\-то конеч\-но-ад\-ди\-тив\-ной} (чистым зарядом, чистым средним), если любая счетно-адди\-тив\-ная мера $\lambda$, удовлетворяющая условию $0\le\lambda\le\mu$ , тождественно равна нулю. 

Знакопеременная мера $\mu$ называется 
{\it чисто конечно-аддитивной}, если в её разложении Жордана $\mu= \mu^{+} - \mu^{-}$ обе неотрицательные меры $\mu^{+}$ и $\mu^{-}$ чисто конечно-аддитивны.
\end{definition}

Если мера $\mu$ чисто конечно-аддитивна, то она равна нулю на каждом одноточечном множестве: $\mu(\{x\}) = 0, \forall x \in X$. Обратное, вообще говоря, неверно (например, для меры Лебега на отрезке $[0, 1]$).

Если счётную аддитивность трактовать как непрерывность меры -- функции множеств, то чисто конечно-аддитивная мера $\mu$ - это ``разрывная мера''. Она существует тогда и только тогда, когда существует последовательность множеств $K_n\in\Sigma, n\in N$, $K_1\supset K_2\supset\dots$, $\lim K_n=\cap_{n=1}^{\infty} K_n=\O$, такая, что $\mu (K_n)\equiv\mu (X)$, т.е. $\lim\mu (K_n)=\mu(X)\ne 0=\mu(\lim K_n)=\mu(\O)$ (это следует из Теоремы 1.22 ~\cite{YoHew1}). \\

\begin{theorem}\label{t1.1} 
~\cite{Ale2}, ~\cite{YoHew1}. (Теорема Александрова-Иосиды-Хьюитта).

Любая конечно-аддитивная мера $\mu$ однозначно разлагается в сумму $\mu =\mu_1 + \mu_2$, где $\mu_1$ - счетно-аддитивная, а $\mu_2$ - чисто конечно-аддитивная меры.
\end{theorem}

Комментарии: 

В статье А. Д. Александрова \cite{Ale2} (глава III, \S 13, теоремы 1, 2, 3, 4, 5, стр. 614 -- 621, и Резюме на русск. яз., стр. 627 -- 628) для некоторых топологических пространств доказывается, что  каждая регулярная конечно-аддитивная мера (А.Д. Александров первый ввёл термин ``заряд'' для таких мер) единственным образом представима в виде суммы счётно-аддитивной меры и ``сингулярной'' конечно-аддитивной меры, равной нулю на каждом компакте. Нетрудно увидеть, что такая ``сингулярная'' мера -- это чисто конечно\--ад\-ди\-тив\-ная мера (для интервала $(0, 1)$ -- это очевидно) по Определению 1.6. 

Спустя 10 лет в фундаментальной работе К. Иосиды и Э. Хьюитта ~\cite{YoHew1} было введено указанное выше понятие чисто конечно\--ад\-ди\-тив\-ной меры и доказано в общем случае существование и единственность разложения любой конечно-аддитивной меры в сумму своей счётно\--ад\-ди\-тив\-ной и чисто конечно\--ад\-ди\-тив\-ной компонент. Однако, ссылок на теоремы (и на статьи) А.Д. Александрова там нет. 

Эта теорема в формулировке ~\cite{YoHew1} получила в литературе название ``Разложение Иосиды-Хьюитта''. В монографии ~\cite[глава III, пункт 15]{DS1} (в блоке ``Конечно\--ад\-ди\-тив\-ные функции множеств'', стр. 255, в оригинале книги ~\cite{DS1} на англ. языке это стр. 233) указываются оба разложения из ~\cite{Ale2} и из ~\cite{YoHew1} с упоминанием фамилий всех их авторов.

Мы полагаем, что есть все основания считать А. Д. Александрова автором первой частной версии обсуждаемого разложения.

Выше мы приводим эту Теорему 1.1 в формулировке ~\cite{YoHew1}, но под названием ``Разложение Александрова-Иосиды-Хьюитта''.
\medskip

 В настоящей статье рассматриваются также банаховы пространства ограниченных мер $\mu \colon \Sigma\to R$ , с нормой, равной полной вариации меры $\mu$ (но можно использовать и эквивалентную sup-норму): $ba(X,\Sigma)$ - пространство конечно-аддитивных мер, $ca(X,\Sigma)$ - пространство счетно-аддитивных мер.  Если $\mu\ge{0}$, то норма $\| \mu \|=\mu(X)$. Меру, тождественно равную нулю, можно формально считать и счётно-аддитивной и чисто конечно-аддитивной. \medskip
 
Чисто конечно-аддитивные меры также образуют банахово пространство $pfa(X,\Sigma)$ с той же нормой, $ba(X,\Sigma) =ca(X,\Sigma) \oplus pfa(X,\Sigma)$. 

Если $X$ - топологическое пространство с топологией $\tau$ ($\tau$ - множество всех открытых подмножеств в $X$), то в качестве $\Sigma$ берём борелевскую сигма-алгебру $\mathfrak B$, порожденную топологией $\tau$.

Для топологических пространств $X$ используется также обозначение $C(X)$ - банахово пространство ограниченных непрерывных функций $f\colon X \to R$ с $sup$-нормой. \medskip

Обозначим множества мер: 

$S_{ba}=\{\mu\in{ba(X,\Sigma)}\colon \mu\ge{0}, ||\mu||=1\},$
 
$S_{ca}=\{\mu\in{ca(X,\Sigma)}\colon \mu\ge{0}, ||\mu||=1\},$

$S_{pfa}=\{\mu\in{pfa(X,\Sigma)}\colon \mu\ge{0}, ||\mu||=1\}.$

Все меры из этих множеств будем называть {\it вероятностными}.\medskip 

\begin{remark}\label{r2}
Мы изначально требуем, чтобы семейство подмножеств $\Sigma$ в исходном произвольном множестве $X$ являлось сигма-алгеброй. Однако, ряд результатов и рассмотрений настоящей статьи останется верным и для случая, когда $\Sigma$ - алгебра подмножеств в $X$. 
\end{remark}

{\bf Примеры}. Приведём два примера чисто конечно-аддитивных мер.

1. Пусть  $X=[0, 1] \subset R$ ($R$ - числовая прямая), $\Sigma = \mathfrak B$. Существует (доказано) такая конечно-аддитивная мера $\mu \colon \mathfrak B \to R$, $\mu \in S_{ba}$, что для любого $\varepsilon > 0$ выполняется: 
$$\mu((0, \varepsilon)) = 1, \mu([\varepsilon, 1]) = 0, \mu(\{0\}) = 0.$$
Неформально можно сказать, что мера $\mu$ фиксирует единичную массу как угодно близко к нулю (справа), но не в нуле. Согласно ~\cite{YoHew1}, такая мера является чисто конечно-аддитивной, но она не единственная. Известно, что мощность семейства таких мер, расположенных ``около нуля (справа)'' не менее, чем $2^{2^{\aleph{0}}} = 2^{c}$ (гиперконтинуум). И такое же семейство чисто конечно-аддитивных мер существует ``около каждой точки $x_{0} \in [0, 1]$ (справа, или слева, или и там, и там)''.
  
2. Пусть  $X = R, \Sigma = \mathfrak B$. Существует (доказано) такая конечно-аддитивная мера $\mu \colon \mathfrak B \to R$, $\mu \in S_{ba}$, что для любого $x \in R$ выполняется:
$$\mu((x, \infty)) = 1, \mu((-\infty, x)) = 0, \mu(\{x\}) = 0.$$
Можно сказать, что мера $\mu$ фиксирует единичную массу как угодно далеко, ``около бесконечности''. Эта мера также чисто конечно-аддитивна. И таких мер также очень много.

Существуют и более сложно устроенные чисто конечно-аддитивные меры, в том числе очень ``похожие'' на меру Лебега.

Подробное изложение основ общей теории конечно-аддитивных мер содержится в монографии K. P. S. Bhaskara Rao, M. Bhaskara Rao ~\cite{RaoRao1}. В ней такие меры называются зарядами (chardge), а разложение мер из Теоремы 1.1 называется ``разложением Иосиды-Хьюитта'' (ссылок на статьи и теоремы А.Д. Александрова в ~\cite{RaoRao1} нет).


\section{Марковские операторы и конечно-аддитивные меры}\label{s2}

В настоящей статье рассматриваются общие цепи Маркова (General Markov Chains), понимаемые как случайный марковский процесс с дискретным временем в произвольном фазовом пространстве и однородный по времени. Заметим, что и сегодня термин ``цепи Маркова'' трактуется в литературе, в том числе в базовых учебниках, весьма различно.

Цепи Маркова (ЦМ) на (фазовом) измеримом пространстве $(X,\Sigma)$ задаются своей переходной функцией (вероятностью) $p(x,E), x\in X, E\in\Sigma$, при обычных условиях: \medskip

1) $0\le p(x,E) \le{1}, p(x,X)=1, \forall{x}\in{X}, \forall{E}\in\Sigma$;

2) $p(\cdot,E)\in{B(X,\Sigma)}, \forall{E}\in\Sigma$;

3) $p(x,\cdot)\in{ca(X,\Sigma)}, \forall{x}\in{X}$. \medskip

Числовое значение функции $p(x, E)$ - это вероятность того, что система перейдёт из точки $x\in X$ в множество $E \in \Sigma$ за один шаг (за единицу времени).

Подчеркнем, что {\itпереходная функция у нас является счетно-аддитивной мерой} по второму аргументу, т.е. мы рассматриваем {\it классические ЦМ.} \medskip

Переходная функция порождает два марковских линейных ограниченных положительных интегральных оператора:\medskip
 
$T: B(X,\Sigma)\to{B(X,\Sigma)}, (Tf)(x)=Tf(x)=\int_X f(y)p(x,dy),$

$\forall{f\in{B(X,\Sigma)}},\forall{x}\in{X};$ \medskip

$A: ca(X,\Sigma)\to{ca(X,\Sigma)}, (A\mu)(E)= A\mu(E)=\int_X p(x,E)\mu(dx),$

$\forall{\mu\in{ca(X,\Sigma)}},\forall{E}\in\Sigma.$ \medskip

Пусть начальная мера $\mu_0\in{S_{ca}}$. Тогда итерационная последовательность счетно-аддитивных вероятностных мер $\mu_{n}=A\mu_{n-1}\in{S_{ca}}, n\in{N}$, обычно и отождествляется с цепью Маркова. Будем называть $\{ \mu_{n}\}$ {\it марковской последовательностью мер}. \medskip

Подобный операторный подход в теории цепей Маркова впервые был сформулирован в работах Н. Крылова и Н. Боголюбова ~\cite{KrBo1}, ~\cite{KrBo2} (1937 год), и детально был разработан в работе К. Иосида и С. Какутани ~\cite{YoKa1} (1941 год), что и отражено в ее названии. Важно подчеркнуть, что в этих трёх работах все используемые меры, в том числе переходная функция $p(x,\cdot)$, счетно-аддитивны. Описываемое ниже продолжение оператора A на пространство конечно-аддитивных мер было произведено рядом авторов намного позже.

Интересно также отметить, что один из соавторов работы ~\cite{YoKa1} К. Иосида, позже (1952 год) стал одним из соавторов работы ~\cite{YoHew1} (совместно с Э. Хьюитт), в которой разработаны основы современной общей теории конечно-аддитивных мер. Однако, цепи Маркова, как возможный объект приложения этой новой теории, в ~\cite{YoHew1} не упонимаются. 

При изучении общих цепей Маркова в рамках операторного подхода обычно в явном виде не используются случайные величины (элементы) $\xi_{n}$, соответствующие марковской последовательности мер $\mu_{n}, n \in N,$ поскольку их значения в произвольном измеримом пространстве $(X,\Sigma)$ нельзя складывать, умножать на число и т.п. Однако, в частных случаях такие марковские случайные элементы (величины) полезны, мы их используем в примерах в \S 6.

В настоящей статье, как и в предыдущих работах автора ~\cite{Zhd01}, ~\cite{Zhd02} и в других, изучаются некоторые проблемы, лежащие в области соприкосновения операторной теории общих цепей Маркова и теории конечно-аддитивных мер. \medskip

 Топологически сопряженным к пространству $B(X,\Sigma)$ является (изоморфно) пространство конечно-аддитивных мер: $B^*(X,\Sigma)=ba(X,\Sigma)$ (см., например, ~\cite{DS1}). При этом топологически сопряженным к оператору Т служит оператор $T^*\colon ba(X,\Sigma)\to{ba(X,\Sigma)}$, однозначно определяемый по известному правилу через интегральные ``скалярные произведения'':          

$$
\langle T^*\mu, f \rangle=\langle \mu, Tf\rangle  \text{ для всех } f\in{B(X,\Sigma)} \text{ и } \mu\in{ba(X,\Sigma)}.
$$ 

Оператор $T^*$  является единственным ограниченным продолжением оператора А на все пространство $ba(X,\Sigma)$ с сохранением его аналитического вида, т.е.  
$$
T^*\mu(E)=\int_X p(x,E)\mu(dx), \forall{\mu\in{ba(X,\Sigma)}}, \forall E\in\Sigma.
$$ 

Интеграл от ограниченной измеримой функции по ограниченной конечно\--ад\-ди\-тив\-ной мере строится по той же схеме, как и интеграл Лебега по произвольной ограниченной счетно-аддитивной мере (см., например, ~\cite{DS1}, ~\cite{YoHew1}). 

Оператор $T^*$ имеет собственное инвариантное подпространство $ca(X,\Sigma)$, т.е. \; \; $T^*(ca(X,\Sigma)) \subset{ca(X,\Sigma)}$, на котором он совпадает с изначальным оператором $A$. Теперь конструкция марковских операторов Т и Т* уже функционально замкнута. Мы будем по-прежнему обозначать оператор Т* как А. 

В такой постановке естественно допустить к рассмотрению и марковские последовательности вероятностных конечно-аддитивных мер: 

$$
\mu_0\in{S_{ba}}, \mu_n = A\mu_{n-1}\in{S_{ba}}, n\in{N},
$$ 
сохраняя счетную аддитивность переходной функции $p(x,\cdot)$ по второму аргументу. Несмотря на это обстоятельство, образ $A\mu$ чисто конечно-аддитивной меры $\mu$  может остаться чисто конечно-аддитивным, т. е., вообще говоря, 
$$A[ba(X,\Sigma)] \not \subset{ca(X,\Sigma)}.$$ 

Допустимо и кардинально изменить постановку задачи - разрешить самой переходной функции $p(x,\cdot)$ быть всего лишь конечно-аддитивной мерой. Такие ЦМ тоже изучаются (см. ~\cite{Rama1}, и также ~\cite[глава~II, \S 5]{Zhd02}) и называются {\it``конечно-аддитивными ЦМ''}, но в настоящей статье не рассматриваются. 

Таким образом, в нашем случае уместна следующая терминология: изучаются {\it счетно-аддитивные} цепи Маркова с операторами, заданными на пространстве {\it конечно-аддитивных мер}.\medskip

Пусть даны произвольное множество $X$ и проивольная сигма-алгебра его подмножеств $\Sigma$, содержащая все одноточечные множества. 
На фазовом пространстве $(X, \Sigma)$ задана произвольная (общая) цепь Маркова с переходной счетно-аддитивной функцией $p(x, E)$ и марковскими операторами $T$ и $T^*=A$.

\begin{definition}\label{de2.1}
Если для некоторой меры $\mu \in S_{ba}$ выполняется $A\mu = \mu$, то будем называть такую меру {\it инвариантной} для оператора $A$ (и для цепи Маркова). Инвариантную счётно-аддитивную меру часто называют также {\it стационарным распределением} ЦМ. В терминах функционального анализа инвариантная мера является неподвижной точкой оператора $A$, или его собственным вектором с собственным числом $\lambda =1$.
\end{definition}

Обозначим множества инвариантных вероятностных мер ЦМ:

$\Delta_{ba}=\{\mu\in{S_{ba}}: \mu=A\mu\}$, 

$\Delta_{ca}=\{\mu\in{S_{ca}}: \mu=A\mu\}$, 

$\Delta_{pfa}=\{\mu\in{S_{pfa}}: \mu=A\mu\}$. \medskip

Обозначим $M_{ba}$ - линейное подпространство инвариантных мер ЦМ в пространстве $ba(X,\Sigma)$.  
Очевидно, что$M_{ba}$ порождается множеством $\Delta_{ba}$: $M_{ba}=Sp\Delta_{ba}$. Будем использовать также обозначения $M_{ca}$ и $M_{pfa}$ с аналогичным смыслом.

Линейной размерностью множества $\Delta_{ba}$ будем называть алгебраическую размерность порождённого им линейного пространства $M_{ba}$, и обозначать её $dim \Delta_{ba} = dim M_{ba}$. Пространство $M_{ba}$ может быть конечномерным,  $dim M_{ba} < \infty$, или бесконечномерным, $dim M_{ba} = \infty$.

Классическая счетно-аддитивная цепь Маркова может иметь инвариантные вероятностные счетно-аддитивные меры, а может и не иметь, т.е. возможно $\Delta_{ca}=\O$ (например, у симметричного блуждания на $Z$). \medskip

З. Шидак ~\cite{Si1}, который один из первых в рамках операторного подхода продолжил марковский оператор на пространства конеч\-но\--ад\-ди\-тив\-ных мер, доказал следующие две важные теоремы.

\begin{theorem}\label{t2.1} 
(\cite[теорема 2.2]{Si1})
Любая счётно-аддитивная цепь Маркова на произвольном измеримом пространстве $(X,\Sigma)$ имеет хотя бы одну инвариантную конечно-аддитивную меру, т.е. всегда $\Delta_{ba}\ne\O$.
\end{theorem}

 Этот результат был затем коротко доказан в работе автора ~\cite{Zhd01}, как простое следствие из теоремы Крейна-Рутмана (\cite[теорема 3.1]{KR1}). 

\begin{theorem}\label{t2.2}
(\cite[теорема 2.5]{Si1}). Если конечно-аддитивная мера $\mu$ для произвольной ЦМ инвариантна $A\mu=\mu$, и $\mu=\mu_1+\mu_2$ - ее разложение на счетно-аддитивную и чисто конечно-аддитивную компоненты, то каждая из них также инвариантна: $A\mu_1=\mu_1$, $A\mu_2=\mu_2$. 
\end{theorem}

Следовательно, $\Delta_{ba}=co\{\Delta_{ca},\Delta_{pfa}\}$, $M_{ba}=M_{ca} \oplus M_{pfa},$ и достаточно изучать инвариантные меры из $\Delta_{ca}$ и из $\Delta_{pfa}$ по отдельности. \\
 
 В работах Ш. Р. Фогеля ~\cite{Fo0}, ~\cite{Fo1} и в других развивается операторный подход для изучения феллеровских цепей Маркова на локально-компактных топологических пространствах. Используется пара марковских операторов  $T\colon C(X) \to C(X)$ и $T^{*} = A\colon rba(X, \mathfrak B) \to rba(X, \mathfrak B)$, где
 $rba(X, \mathfrak B)$ обозначает пространство регулярных ограниченных конечно-аддитивных мер (зарядов) на $(X, \mathfrak B)$, причём  $C^{*}(X)=rba(X, \mathfrak B)$. В ~\cite{Fo1}, в частности, доказывается (для положительных операторов), что если $\mu=A\mu \in rba(X, \mathfrak B)$ и 
 $\mu =\mu_1 + \mu_2$ -- разложение меры $\mu$ две компоненты: $\mu_1$ -- регулярную счётно-аддитивную и $\mu_2$ -- регулярную чисто конечно-аддитивную, то $A\mu_{1}=\mu_{1}$ и $A\mu_{2}=\mu_{2}$.
 
 В книге Ш. Р. Фогеля ~\cite{Fo2} (и в других его статьях) рассматриваются также цепи Маркова и марковские операторы $T$ и $T^{*}=A$ на пространствах  $L_{1}(X, \Sigma, m)$ и $L_{1}^{*}=L_{\infty}(X, \Sigma, m)$, где $m$ -- заранее заданная фиксированная $\sigma$-конечная счётно-аддитивная мера. В этой книге также используются конечно-аддитивные меры. 
 
 В настоящей статье мы рассматриваем общие цепи Маркова и не выделяем отдельно случай топологического фазового пространства $(X, \mathfrak B)$.
 Однако, представленные далее наши результаты при некоторой корректировке будут верны и для таких пространств.
 
  Мы также не предполагаем на фазовом пространстве $(X, \Sigma)$ никакой априорно заданной фиксированной меры $m$, и не используем пространства  $L_{1}, L_{\infty}$, а также  $L_{p}, 1<p<\infty$. Отметим, что кроме Ш. Р. Фогеля и его учеников, ещё ряд не названных здесь авторов работает с цепями Маркова с использованием пространств  $L_{1}, L_{\infty}$ и $L_{p}, 1<p<\infty$.
  

\section{Порядковые свойства пространств инвариантных мер}\label{s3}

Приведем некоторые сведения из теории полуупорядоченных линейных пространств (решёток), применительно к нашим пространствам мер, ориентируясь на источники ~\cite{KR1}, ~\cite{Vu1}, ~\cite{Bir1}, ~\cite{YoHew1}, ~\cite{RaoRao1}.

\begin{definition}\label{de3.1}
Пусть $L$ - линейное (векторное) нормированное пространство. Его замкнутое подмножество $G \subset L, G \ne \O$ называется {\it конусом} положительных (неотрицательных) элементов в $L$, если оно обладает следующими свойствами:

1) если $x \in G$, то  $\lambda x \in G$ при всех $\lambda \ge 0$,

2) если $x, y \in G$, то  $x+y \in G$,

3) если $x \in G$, то  $-x \ne G$.
\end{definition}

Каждый конус  $G \subset L$ порождает частичное упорядочивание (полуупорядочивание) в  $L$  для  $x_{1}, x_{2} \in L$:  $x \le y$, если  $y-x \in K$.

Пусть $L$ - вещественное линейное (векторное) нормированное пространство, полуупорядоченное по конусу $G$ своих положительных элементов.  

\begin{definition}\label{de3.2}
Конус $G$ называется {\it воспроизводящим}, если  $G \ominus G = L$, т. е. любой элемент из $L$ представим в виде разности двух элементов из $G$.
\end{definition}

\begin{definition}\label{de3.3}
Конус $G$ называется {\it миниэдральным}, если для любых $x, y \in G$ существуют  $z = sup(x, y) = x \vee y \in G, \; \tilde{z} = inf(x, y) = x \wedge y \in G$.
\end{definition}
  
  Если конус $G$  миниэдральный и воспроизводящий, то супремум и инфимум существует и для любых  $x, y \in L$.
  
\begin{definition}\label{de3.4}
 ~\cite{Vu1} Линейное пространство  $L$, полуупорядоченное  по миниэдральному воспроизводящему конусу  $G$, называется  {\it $K$-линеалом}. В терминах ~\cite{Bir1} это {\it векторная решётка}.
\end{definition}
  
  Любой элемент $x$ $K$-линеала $L$ можно однозначно разложить в сумму $x = x_{+} - x_{-}$ своих положительной и отрицательной частей:  $x_{+} = x \vee 0,  \; x_{-} = (- x) \vee 0$ . Элемент $|x| = x_{+} + x_{-}$ называется  {\it модулем} элемента  $x$.

Пространство мер $ba(X,\Sigma)$ полуупорядочено по естественному отношению порядка:

для  $\mu_{1}, \mu_{2} \in ba(X,\Sigma)$ пишем  $\mu_{1} \le \mu_{2}$, если $\mu_{1}(E) \le \mu_{2}(E)$ для всех $E \in \Sigma$.
Это отношение порядка используется также и для подпространств мер $ca(X,\Sigma)$ и  $pfa(X,\Sigma)$.

Выделим в этих пространствах конусы всех положительных  элементов  $G_{ba}$, $G_{ca}$, $G_{pfa}$. Все эти конусы являются воспроизводящими, т. е.  порождают разностями все пространство. 

В ~\cite{YoHew1} приводится теорема 1.11, принадлежащая С. Бохнеру и Р. С. Филлипсу (1941) о том, что естественному порядку в пространстве $ba(X, \Sigma)$ (а значит и в подпространствах $ca(X, \Sigma)$ и $pfa(X, \Sigma)$) соответствует следующая конструкция:
$$
(\mu_{1} \land \mu_{2})(E) = inf \{ \mu_{1}, \mu_{2} \}(E) = inf ( \mu_{1}(C) + \mu_{2}(E \cap C^{\prime})), \text{ где } C^{\prime}=X\setminus C, 
$$ 
и инфимум берется по всем множествам  $C \subset E, C \in \Sigma$;
$$
(\mu_{1} \lor \mu_{2}) =sup \{ \mu_{1}, \mu_{2} \} = -((-\mu_{1}) \land (-\mu_{2})).  
$$

При этом $\mu_{1} \land \mu_{2}, \mu_{1} \lor \mu_{2} \in ba(X,\Sigma)$.

Таким образом, с этой конструкцией пространство $ba(X, \Sigma)$ является $K$-ли\-неа\-лом (решёткой) и конус $G_{ba}$ является миниэдральным. 

В ~\cite{YoHew1}, в теореме 1.14, доказывается, что, если меры $\mu_{1}$ и $\mu_{2}$ счётно-аддитивны, то и меры $\mu_{1} \land \mu_{2}, \mu_{1} \lor \mu_{2}$ счётно-аддитивны. Там же, теорема 1.17, доказывается, что и для чисто конечно-аддитивных мер $\mu_{1}$ и $\mu_{2}$ их инфимум и супремум также чисто конечно-аддитивны. Таким образом,  $K$-линеалами (решётками) являются и подпространства $ca(X, \Sigma)$ и $pfa(X, \Sigma)$, а их конусы положительных элементов -- миниэдральны.\medskip

Нам понадобится ещё несколько определений.

\begin{definition}\label{de3.5}
 (\cite[глава  III, \S 10, определение  III.10.1, стр. 81]{Vu1})
$K$-линеал $X$
 называется {\it архимедовым}, если в нем выполнен принцип Архимеда: из того, что для некоторого $x \ge 0 \; (x \in X)$ множество  $n \cdot x$, $n \in N,$ всех его ``кратных'' ограничено, вытекает равенство  $x=0$.
\end{definition} 

 \begin{lemma}\label{l3.1}
 $K$ - линеал $ba(X, \Sigma)$  и любое его линейное подпространство (конечномерное или бесконечномерное) является архимедовым.
\end{lemma}
 
 Доказательство легко следует из того, что $ba(X, \Sigma)$ линейное нормированное пространство и все меры $\mu \in ba(X, \Sigma)$ ограничены.
 
 \begin{corollary}\label{cor3.1}
 $K$ - линеалы $ca(X, \Sigma)$ и $pfa(X, \Sigma)$ и любые их линейные подпространства (конечномерные или бесконечномерные) являются архимедовыми.
\end{corollary}

\begin{definition}\label{de3.6}
Два неотрицательных элемента $x_{1}$ и $x_{2}$ из $K$ - линеала $X$ называются {\it дизъюнктными}, если  $x_{1} \land x_{2} =0$, что обозначается  $x_{1}dx_{2}$. Этот термин мы применяем и для соответствующих пар мер $\mu_{1}, \mu_{2} \in ba(X, \Sigma)$
\end{definition}

\begin{definition}\label{de3.7}
 (см., например, ~\cite[глава III, п. 4.12, стр. 146]{DS1}). Неотрицательные меры $\mu_{1}$ и $\mu_{2}$  называются {\it сингулярными}, если существуют множества $D_{1}, D_{2} \in \Sigma$, такие, что $\mu_{1}(D_{1}) = \mu_{1}(X), \mu_{2}(D_{2}) = \mu_{2}(X)$ и  $D_{1} \cap D_{2} = \O$, что обозначается $\mu_{1}\perp \mu_{2}$. При этом множества $D_{1}$ и $D_{2}$ называют иногда {\it носителями} мер $\mu_{1}$ и $\mu_{2}$ или множествами полной меры для $\mu_{1}$ и $\mu_{2}$ соответственно. Определение носителя меры здесь неоднозначно (более точное определение носителя меры возможно только в топологическом пространстве  $X$, что нам сейчас не нужно, поскольку множество $X$ произвольно).
\end{definition}
 
 Нетрудно проверить, что если меры $\mu_{1}$ и $\mu_{2}$  сингулярны, то они и дизъюнктны, т. е. $\mu_{1}\perp \mu_{2} \Rightarrow \mu_{1}d\mu_{2}$. Однако, обратное, вообще говоря, неверно (в ~\cite{RaoRao1} приведены соответствующие примеры).
 
\begin{remark}\label{r4} 
В книге ~\cite{RaoRao1} (пункты 6.1.14 и 6.1.15, стр. 164) наши дизъюнктные пары мер (заряды) называются {\it сингулярными}, а наши сингулярные меры называются {\it сильно сингулярными}. В статье ~\cite{YoHew1} меры, для которых $\mu_{1} \wedge \mu_{2}=0$, вообще никаким специальным словом не называются, а сингулярность понимается так же, как и у нас (точнее - у нас так же, как у них).
\end{remark}
 
 Нам понадобится также следующая Теорема 1.21 из ~\cite{YoHew1}, которую мы приведем в нашей символике.
 
\begin{theorem}\label{t3.1}
 ~\cite{YoHew1}.
 
  1. Пусть две конечно-аддитивные меры $\mu_{1},\mu_{2} \in S_{ba}$ дизъюнктны. Тогда для любых чисел $\alpha, \beta >0$ существует множество $E\in \Sigma$ такое, что  $\mu_{1}(E) \le \alpha$ и $\mu_{2}(X \setminus E) \le \beta$.  
 Здесь можно положить $\varepsilon = \alpha = \beta$ и считать $\varepsilon$ ''как угодно малым''. 
  
 2. Если хотя бы одна из двух конечно-аддитивных мер $\mu_{1}$ и $\mu_{2}$ не является счетно-аддитивной, то из их дизъюнктности не следует сингулярность.
 
  3. Счетно-аддитивные меры $\mu_{1}, \mu_{2}$ дизъюнктны тогда и только тогда, когда они сингулярны.\medskip
\end{theorem}

В теореме 1.16 из ~\cite{YoHew1} утверждается также, что каждая чисто конечно-аддитив\-ная мера дизъюнктна с любой счётно-аддитивной мерой.  
 
 \begin{definition}\label{de3.8}
(\cite[глава  III, \S 13, определение  III.13.1, стр. 88]{Vu1}).
 Элемент $x$ $K$-линеала $X$ называется {\it дискретным}, если не существует дизъюнктных между собой элементов y>0 и z>0, таких, что $y \le |x|$ и  $z \le |x|$. Это означает, что если y и z таковы, что  $0 \le y \le |x|$, $0 \le z \le |x|$ и  $ydz$, то или $y=0$ или $z=0$ (последние два неравенства можно заменить на равенство $y + z = |x|$). \medskip
\end{definition}
 
 Теперь приведем одну классическую структурную теорему А. И. Юдина (1939 год), представленную вместе с нетривиальным доказательством в 
(\cite[глава  III, \S 14, определение  III.14.1, стр. 89-91]{Vu1}), которую мы далее будем использовать. Теорема А. И. Юдина цитируется во многих работах по теории решёток, в том числе в ~\cite{Bir1}.

\begin{theorem}\label{t3.2} 
~\cite{Vu1}).
Любой  $n$-мерный архимедов  $K$-линеал $X$ (архимедова векторная решётка) алгебраически и структурно изоморфен $K$-линеалу  $R^{n}$.
\end{theorem}

Здесь $R^{n}$ - обычное $n$-мерное евклидовое  пространство с ортонормированным базисом векторов  $H=\{e_{1}, e_{2}, ..., e_{n}\}, e_{1} = (1, 0, ..., 0), e_{2} = (0, 1, ..., 0), ..., e_{n} = (0, 0, ..., 1)$.

В доказательстве этой теоремы в ~\cite{Vu1} конкретизируется и уточняется ее утверждение. 
Эту конкретизацию мы используем ниже при применении теоремы А. И. Юдина, а здесь - опустим. \medskip

 Для дальнейшего исследования нам нужно уточнить порядковые свойства пространств инвариантных мер для общих цепей Маркова  $M_{ba}$, $M_{ca}$  и $M_{pfa}$. Они являются линейными подпространствами пространств мер $ba(X, \Sigma)$,  $ca(X, \Sigma)$ и $pfa(X, \Sigma)$, которые, как указано выше, являются $K$-линеалами (векторными решётками). Возникает естественный вопрос -- а являются ли $M_{ba}$, $M_{ca}$  и $M_{pfa}$ и подрешётками в этих пространствах?  Уточняем вопрос.

Для любых $\mu_1 , \mu_2 \in M_{ba}$  (или  $M_{ca}$, или $M_{pfa}$), 
как элементов решетки $ba(X, \Sigma)$ $(\text{или } ca(X, \Sigma), \text{или } pfa(X, \Sigma))$,
существует $sup(\mu_1 , \mu_2)$ и $inf(\mu_1 , \mu_2)$, принадлежащие
$ba(X, \Sigma)$ $(\text{или } ca(X, \Sigma), \text{или } pfa(X, \Sigma))$.
А будут ли эти супремум или инфимум инвариантными для оператора
$A$ цепи Маркова, т. е. попадают ли они
в подпространства $M_{ba}$ (или $M_{ca}$, или $M_{pfa}$)?
Это не такой уж простой вопрос, как может показаться на первый взгляд.

Некоторый ответ на этот вопрос даётся в одной теореме Г. Биркгофа 
(\cite[глава XVI, пункт 7, теорема 12, стр. 502]{Bir1}), где говорится (и доказывается), что множество F точек произвольного
$(L)$-пространства, неподвижных относительно какого-нибудь оператора перехода $A$ (у автора $T$), метрически замкнуто, оно является подпространством и подрешеткой. Не будем давать определение $(L)$-пространства (см. ~\cite[глава XV, пункт 15, стр. 479]{Bir1}), заметив только, что наши пространства мер $M_{ba}$, $M_{ca}$ и $M_{pfa}$ таковыми являются. Под оператором перехода понимается ``линейный оператор $A$, переводящий вероятностные распределения в вероятностные''. Из контекста видно, что ``вероятностное распределение'' - это ``вероятностная мера'', которая по определению 
(\cite[глава XI, пункт 5, стр. 339]{Bir1}) является счетно-аддитивной.

Таким образом, данная теорема Биркгофа применима лишь для случая марковского оператора $A\colon ca(X, \Sigma) \to ca(X, \Sigma)$ (оператора перехода), заданного на пространстве счётно-аддитивных мер. В этом случае подпространство $M_{ca}$ является подрешёткой. 
Если же мы продолжим марковский оператор $A$ до оператора $T^{*}=A\colon ba(X, \Sigma) \to ba(X, \Sigma)$  на пространстве конечно-аддитивных мер, то здесь данная теорема Биркгофа (формально) уже не применима.
Однако, то,  что подпространство инвариантных  конечно-аддитивных мер $M_{ba}$  является  и подрешёткой в $ba(X, \Sigma)$ доказывается, в частности, (без связи с теоремой Биркгофа) в нашей статье (\cite[глава ~II, \S 6, теорема 6.4, стр. 74]{Zhd01}). 

\begin{theorem}\label{t3.3} 
Пусть произвольная цепь Маркова имеет  конечномерное пространство $M_{ba}$ всех своих инвариантных конечно аддитивных мер,  т. е. $dim M_{ba} = dim \Delta_{ba} = m < \infty$.  Тогда пространство $M_{ba}$ алгебраически и структурно (решёточно)  изоморфно евклидовому  пространству $R^m$ и имеет линейно независимый базис инвариантных мер $H_{ba}=\{ \mu_1 ,  \mu_2 , ... , \mu_m \} $, такой, что $\mu_i \in S_{ba}$, где все меры $\mu_i$ дискретны и являются либо счётно-аддитивными, либо чисто конечно-аддитивными при $i = 1, 2, ..., m$, они попарно дизъюнктны, т. е. $\mu_i \wedge \mu_j = 0, \mu_i d \mu_j $ при $i \neq j$.

При этом любая инвариантная мера $\mu \in M_{ba}$ единственным образом представима в виде $$\mu = \sum_{i=1}^{m} \alpha_i \mu_i ,$$ где $\mu_i \in H_{ba}, \alpha_i \in R$ при $i=1, 2, ..., m$.
\end{theorem}

{\bf{Доказательство}}.
Как отмечено выше, в ~\cite[глава II, \S 6, теорема 6.4]{Zhd01} доказано, что пространство $M_{ba}$  конечно-аддитивных инвариантных мер для ЦМ является $K$-линеалом (векторной решёткой). По лемме 3.1. пространство $M_{ba}$ является архимедовым. По  условию теоремы пространство $M_{ba}$ конечномерно.  Следовательно,  применима теорема Юдина,  представленная выше как теорема 3.2 ~\cite{Vu1}.
Заменяем символику в общей теореме 3.2  Юдина и  в комментариях к ней на символику  нашего частного случая для пространства мер $M_{ba}$. 

Остаётся только пояснить дискретность базисных инвариантных мер $\mu_i, i=1, 2, ..., m$. Пусть базисная конечно-аддитивная мера $\mu \in H_{ba}$ и $\mu = \mu_{ca} + \mu_{pfa}$ её  разложение Александрова-Иосиды-Хьюитта на счётно-аддитивную  и чисто конечно-аддитивную компоненты. Тогда по теореме 
Шидака (см. \cite[теорема 2.5]{Si1}) обе меры $\mu_{ca}$ и $\mu_{pfa}$ также инвариантны. Тогда меры $\mu_{ca}$ и $\mu_{pfa}$ дизъюнктны (см. ~\cite[теорема 1.16]{YoHew1}), и $\mu_{ca} \le \mu$, $\mu_{pfa} \le \mu$. Базисные меры $\mu_i$, в том числе $\mu$, по теореме 3.2  Юдина дискретны, что по определению 3.8 означает, что или $\mu_{ca}=0$, или $\mu_{pfa}=0$. Теорема доказана.\medskip

Базис $H_{ba}$ в  пространстве $M_{ba}$,  существующий по утверждению Теоремы  3.3,  будем называть {\it дизъюнктным базисом}.
Аналог Теоремы 3.3  для пространства $M_{ca}$  инвариантных счётно-аддитивных мер  был получен в нашей работе (\cite[глава II, \S 6, теорема 6.3, стр. 73]{Zhd01}).  Эта теорема также опирается на теорему А. И. Юдина. Там счётная аддитивность инвариантных мер позволяет из дизъюнктности базисных мер получить их сингулярность.

Одна из основных решаемых задач в следующих параграфах - это доказать, хотя бы и с оговорками, что все инвариантные меры ЦМ при $dimM_{ba} < \infty$ являются счётно-аддитивными.

 \section{Теоремы о составе множества инвариантных мер и квазикомпактности цепи Маркова}\label{s4}

В работах автора ~\cite{Zhd01} и ~\cite{Zhd02} было показано, что конечная размерность и состав множества {инвариантных} {конечно-адди\-тив\-ных мер} марковского оператора А тесным образом связаны с одним из центральных вопросов эргодической теории ЦМ, а именно, с известными условиями квазикомпактности марковских операторов (см., например, ~\cite{YoKa1} и ~\cite[глава V, \S 5]{Doob1}).\medskip

\begin{definition}\label{de4.1}
Марковский оператор $A\colon ca(X, \Sigma) \to ca(X, \Sigma)$ называется  {\it квазикомпактным} (квази-вполне непрерывным), если выполнено условие \\
$
(K)
\begin{cases}
\text{{ существуют {\it компактный} (вполне непрерывный) оператор $A_1$}} \\
\text{{ (переводящий каждое ограниченное множество в предкомпактное)  }}\\
\text{{ и целое число $k\ge1$ такие, что $||A^k-A_1||<1$.  }} 
 \end{cases}
$ 

В этом случае и саму ЦМ будем называть {\it квазикомпактной}. \medskip
\end{definition}

 В ~\cite{YoKa1} доказывается, что если оператор $A$ квазикомпактен, то цепь Маркова имеет конечное число сингулярных инвариантных счётно-аддитивных вероятностных мер и эргодические средние по Чезаро от марковской последовательности мер равномерно по начальным счётно-аддитивным мерам и экспоненциально быстро сходятся в некотором смысле к ним (к соответствующему проектору) в метрической топологии. Это означает, что верна равномерная (средняя) эргодическая теорема.

Отметим, что в наших терминах для квазикомпактных ЦМ $\Delta_{ca}\ne \O$ и $dim\Delta_{ca}=n<\infty$. Однако, обратное, вообще говоря, не верно.
 
 Если ЦМ не является квазикомпактной, то возможны все три случая: $\Delta_{ca}=\O$, $dim\Delta_{ca}=n<\infty$ и $dim\Delta_{ca}=\infty$. \medskip
 
 Справка.
 
Условие квазикомпактности для марковских операторов $(K)$ ввели в 1937 году Н.М. Крылов и Н.Н. Боголюбов в своих двух кратких заметках, уже упомянутых выше, ~\cite{KrBo1} и ~\cite{KrBo2} (точнее сказать - в двух частях одного сообщения). В этом сообщении формулируются условия квазикомпактности для марковского оператора $T\colon B(X,\Sigma) \to B(X,\Sigma)$ и очень кратко показывается, что этого условия достаточно для того, чтобы средние по Чезаро от последовательности функций $f_{n}=Tf_{n-1}=T^{n}f_{0}$ сильно и равномерно по начальным функциям $f_{0}$ сходились к проектору ($T_{\infty}$) на конечномерное пространство инвариантных (гармонических) и сингулярных функций оператора $T$. Используются также инвариантные счётно-аддитивные меры для двойственного марковского оператора $A\colon ca(X,\Sigma) \to ca(X,\Sigma)$ (без введения символа $A$ или другого для такого оператора). 
В заметках ~\cite{KrBo1}, ~\cite{KrBo2} ещё не используется термин ``квазикомпактность'' и нет маркера $(K)$ для соответствующих условий. Авторы подчеркивают, что фазовое пространство рассматриваемой цепи Маркова $(X,\Sigma)$ максимально общее, а цепь Маркова задаётся произвольной переходной функцией $p(x, E)$, счётно-аддитивной по второму аргументу.
 
Случайные величины, соответствующие марковской последовательности мер, в ~\cite{KrBo1}, ~\cite{KrBo2} упоминаются, но не используются (они и не могут использоваться, так как фазовое пространство $(X,\Sigma)$ -- общее). \medskip

В том же 1937 году была опубликована большая работа В. Деблина в двух частях ~\cite{Doe1} и ~\cite{Doe2}. В этой работе (\cite[глава 2]{Doe2}) для общих цепей Маркова сформулировано условие $(D)$, при выполнении которых цепь Маркова обладает максимальным набором эргодических свойств. Соответствующие утверждения подробно доказываются в ~\cite{Doe2}.

Приведём это условие Деблина $(D)$: 

$
(D)
\begin{cases}
\text{{\it существуют ограниченная мера $\varphi\in{ca(X,\Sigma)}$, $\varphi\ge{0},\varepsilon>0$ и $k\in Z$, $k\ge1$,}}\\
\text{{\it  такие, что из $\varphi(E)\le\varepsilon, E\in\Sigma$, следует $p^k(x,E)\le{1-\varepsilon}$ для всех $x\in{X}$.}}
\end{cases}
$ \medskip

Замечание: верхний индекс $k$ в $p^k$ означает порядок интегральной свертки (итерации) переходной функции, а не её степень. Назовём число $k$ {\it параметром условия} $(D)$.\medskip

В работе К. Иосиды и С. Какутани (~\cite{YoKa1}, 1941 год) детально анализируются и сравниваются эти два условия, им присваиваются маркеры $(K)$ и $(D)$ (в оригиналах их нет), а операторы, удовлетворяющие условию $(K)$, называются ``квази-вполне непрерывными'' (позже их стали называть ``квазикомпактными''). Там же доказывается, что выполнение условия $(D)$ влечёт выполнение условия $(K)$, а условие $(K)$ ``почти'' влечёт условие $(D)$.  

В ~\cite{YoKa1} отмечается, что условие $(K)$ может быть использовано не только для марковского оператора $T$ на пространстве функций, как это сделано в ~\cite{KrBo1}, ~\cite{KrBo2}, но и для ``двойственного'' марковского оператора $A$ на пространстве мер (это мы и повторили выше в определении 4.1), а также для произвольных линейных операторов в банаховых пространствах.  

Обращаем внимание на то, что оператор $A$ в определении 4.1 мы оставляем в классическом варианте условия $(K)$, действующим на пространстве счётно-аддитивных мер $ca(X, \Sigma)$. При других вариантах задания оператора $A$ могут нарушиться некоторые уже имеющиеся импликации с другими условиями. Но это не мешает одновременно рассматривать и продолжение оператора $A=T^{*}$ на пространство $ba(X, \Sigma)$.
 
В работе ~\cite{YoKa1} констатируется, что указанные работы Н.М. Крылова и Н.Н. Боголюбова ~\cite{KrBo1}, ~\cite{KrBo2} заложили основы нового ``теоретико-операторного метода''  в дальнейшем изучении общих цепей Маркова, а работы В. Деблина ~\cite{Doe1},~\cite{Doe2} дали новый ``прямой теоретико-множественный метод'' в этой теории. В настоящей статье используются оба указанных метода. 

В литературе условие $(D)$ нередко формулируют и в других формах (не всегда эквивалентных), в том числе в книге Дж. Дуба ~\cite{Doob1}, в работе ~\cite{YoKa1}, а также в книге Ж. Невё (\cite[глава V, пункт 5.5]{New1}). В нашей работе ~\cite{Zhd02} также дается близкое условие $(\tilde D)$, являющееся некоторой модификацией условия $(D)$. Приведём его. \medskip

Пусть на $(X,\Sigma)$ задана произвольная ЦМ с переходной функцией $p(x,E)$ и марковскими операторами $T$ и $A$. Для любого $m\in N$ определим новую ЦМ с переходной функцией $q_m(x,E)$ и марковскими операторами $T_m$ и $A_m$ по правилам построения средних по Чезаро: 
$$q_m(x,E)=\frac{1}{m}\sum\limits_{k=1}^{m}p^k(x,E),\qquad T_m=\frac{1}{m}\sum\limits_{k=1}^{m}T^k,\qquad A_m=\frac{1}{m}\sum\limits_{k=1}^{m}A^k.$$

Такую ЦМ назовем {\it конечно-осреднённой ЦМ} (по исходной ЦМ).

Сформулируем новое условие квазикомпактности $(\tilde D)$: \medskip \\
$
(\tilde D)
\begin{cases}
\text{{\it существуют ограниченная мера $\varphi\in ca(X,\Sigma)$, $\varphi\ge 0$, $\varepsilon>0$ и $m \in Z, m\ge 1$,  }}\\
\text{{\it такие, что из $\varphi(E)<\varepsilon$, $E\in\Sigma$, следует $q_m(x,E)\le{1-\varepsilon}$ для всех $x\in X$.}}
\end{cases}
$ \medskip

Очевидно, что $(\tilde D)$ является условием Деблина $(D)$ для конечно-осредненнной ЦМ (при фиксированном $m\ge 1$) с параметром $k=1$. Следовательно, если выполнено условие $(\tilde D)$, то оператор $A_m$ является квазикомпактным.

В работе ~\cite[теорема 12.1]{Zhd02} показывается, что если выполнено условие $(D)$, то выполнено и условие $(\tilde D)$. Это означает, что если исходная ЦМ квазикомпактна, то квазикомпактны и все её конечно-осредненные ЦМ. 

Однако, обратное утверждение в работе ~\cite{Zhd02} доказано не было, поэтому и пришлось ввести это условие $(\tilde D)$. Мы смогли получить обращение только теперь в настоящей статье - будет показано ниже.

Обозначим $\tilde\Delta_m$ семейство всех нормированных положительных конечно\--ад\-ди\-тив\-ных инвариантных мер для конечно-осредненной ЦМ с параметром $m$. Очевидно, что $\Delta_{ba}\subset\tilde\Delta_m$ при всех $m\in N$, и, возможно, $\Delta_{ba}\ne\tilde\Delta_m$ при наличии циклических подклассов. Напомним, что для любой ЦМ $\Delta_{ba}\ne \O$. Следовательно, и $\tilde\Delta_m \ne \O$. \medskip

В работе автора ~\cite[теорема 12.2]{Zhd02} было доказано следующее утверждение:\medskip

\begin{theorem}\label{t4.1} 
~\cite{Zhd02}.
Для произвольной ЦМ условие $(\tilde D)$ эквивалентно условию (*):
$$(*) \qquad \Delta_{ba}\subset{ca(X,\Sigma)},$$
которое означает, что все инвариантные конечно-аддитивные меры исходной ЦМ являются счетно-аддитивными, или, другими словами, исходная ЦМ не имеет инвариантных чисто конечно-аддитивных мер.
\end{theorem}

В настоящей работе мы докажем новое утверждение о том, что классическое условие $(D)$, также эквивалентно условию  $(*)$. Но для этого нам понадобятся некоторые сведения из функционального анализа (см., например, ~\cite{DS1}). 

Пространство мер $ba(X, \Sigma)$ является топологически сопряжённым к пространству функций $B(X, \Sigma)$, то есть  $B^{*}(X, \Sigma) = ba(X, \Sigma)$ (с точностью до изоморфизма), что мы уже отмечали. Следовательно, в банаховом пространстве $ba(X, \Sigma)$ можно рассматривать не только сильную (метрическую) топологию, но и $*$-слабую топологию $\tau_{B}$, порожденную предсопряженным пространством $B(X, \Sigma)$. Эта топология задаётся тихоновской базой окрестностей точки (меры) $\eta\in ba(X, \Sigma)$ вида 
$
V(\eta, f_{1}, f_{2}, \ldots, f_{k}, \varepsilon)=\{ \mu \in ba(X, \Sigma) \colon
|\langle f_{i},\mu \rangle-\langle f_{i},\eta \rangle|< \varepsilon,  i = 1, 2, \ldots , k; k\in N, \varepsilon >0 \}
$,  
где $\varepsilon$ и $k$ произвольны, $f_{1}, f_{2}, \ldots, f_{k} \in B(X, \Sigma)$.
Запись $\langle f_{i},\mu \rangle$ обозначает значение функции $f_{i}$, как линейного функционала на мере $\mu$, вычисляемое по формуле  
$$\langle f_{i},\mu \rangle = \int_{X} f_{i}(x)\mu (dx), i = 1, 2, \ldots , k.$$

Если дана некоторая последовательность мер $\{\lambda_{n}\}=\{\lambda_{1}, \lambda_{2}, \ldots \} \subset S_{ba}$, то символом $\mathfrak{M} {\{\lambda_{n}\}}$ обозначим множество всех предельных в $\tau_{B}$-топологии мер последовательности $\{\lambda_{n}\}$. Отметим, что из того, что мера $\eta$ является предельной точкой последовательности мер ${\{\lambda_{n}\}}$, вообще говоря, не следует, что в ${\{\lambda_{n}\}}$ существует подпоследовательность, сходящаяся к $\eta$ в $\tau_{B}$-топологии.\medskip

Пусть на $(X,\Sigma)$ задана ЦМ с оператором $A$. Введем обозначение средних по Чезаро для некоторой начальной меры $\mu \in S_{ba}$:
$$
\lambda_{n} = \lambda^{\mu}_{n} = \frac{1}{n} \sum_{k=1}^{n}A^{k}\mu, n\in N.
$$

Доказательство ряда дальнейших теорем опирается на теорему 7.2 из ~\cite{Zhd01}, доказательство которой в ~\cite{Zhd01} не приводится. Там сказано лишь, что его можно провести по аналогии  с доказательством другой теоремы для феллеровской ЦМ на топологическом пространстве. Однако различие между общими ЦМ и топологическими ЦМ весьма велико, и мы восстанавливаем ниже неопубликованное доказательство теоремы 7.2 из ~\cite{Zhd01}.\medskip
 
\begin{theorem}\label{t4.2}
 Пусть на произвольном $(X,\Sigma)$ заданы произвольная ЦМ и начальная конечно-аддитивная мера $\mu \in ba(X,\Sigma)$, $\mu \in S_{ba}$. Тогда каждая $\tau_{B}$-предельная точка (в топологии $\tau_{B}$) последовательности $\{\lambda_{n}^{\mu}\}$ в $ba(X,\Sigma)$ будет неподвижной точкой оператора $A$,  т.е.  $\mathfrak{M}{\{\lambda_{n}\}} \subset \Delta_{ba}$,  множество таких мер непусто, т.е.  $\mathfrak{M}{\{\lambda_{n}\}} \ne \O$, и  $\mathfrak{M}{\{\lambda_{n}\}}$ -- $\tau_{B}$-компактно. \medskip
\end{theorem}

{\bf{Доказательство.}} 
Выберем некоторую $\mu \in S_{ba}$. Очевидно, что для нее  $\|\lambda_{n}\|=\lambda_{n}(X )=1,\; 
n=1,2,\ldots $, т.е. множество  $\{\lambda_{n}\}$  метрически ограничено в $ba(X, \Sigma)$. Следовательно, 
$\tau_{B}$-замыкание $\{\lambda_{n}\}$ компактно в $\tau_{B}$-топологии (см. ~\cite[глава V, п. 4, следствие 3]{DS1}). 
Согласно Келли (\cite[глава V, теорема 5]{Kel}), любая подпоследовательность  $\lambda_{n_{i}},$  в том 
числе и $\lambda_{n}$, в компакте имеет $\omega$-предельную точку $\eta=\eta\{\lambda_{n_{i}}\}$, такую, что в любой ее окрестности содержится бесконечно много элементов последовательности. Это означает, 
что для каждой $\tau_{B}$-окрестности $V(\eta, f_{1}, f_{2}, \ldots, f_{k}, \varepsilon)$, где $\varepsilon>0$, и $f_{i} \in B(X, \Sigma), i = 1, 2, \ldots , k,$ множество  $\{i: \lambda_{n_{i}} \in V(\eta, f_{1}, f_{2}, \ldots, f_{k},\varepsilon)\}$ бесконечно,  т.е. существует подпоследовательность  
$\{\lambda_{n_{i_{j}}}\}$, \; $\lambda_{n_{i_{j}}} \in V(\eta, f_{1}, f_{2}, \ldots, f_{k}, \varepsilon)$, 
$j=1,2,\ldots$. 

Пусть мера $\eta$ -- $\tau_{B}$-предельна для $\{\lambda_{n}\}$. 

Сделаем следующее преобразование 
$$
A\lambda_{n}=\frac{1}{n}\sum_{k=1}^{n}A^{k+1}\mu+\frac{1}{n}[A^{n+1}\mu-A\mu]=
\lambda_{n}+\frac{1}{n}[A^{n+1}\mu - A\mu].
$$

Пусть  $f \in B(X, \Sigma)$.  Тогда для любого $\varepsilon >0$  существует строго возрастающая последовательность чисел 
$\{n_{i}\}=\{n_{i}\}(f, \varepsilon)$, такая,  что $\lambda_{n_{i}} \in V(\eta, f, Tf, \varepsilon),\; i=1,2,\ldots$.  

Теперь  
$$
|f(\eta)-f(A\eta)|=|f(\eta)-Tf(\eta)| \le 
$$
$$
\le |f(\eta)-f(\lambda_{n_{i}})| + |f(\lambda_{n_{i}}) -Tf(\lambda_{n_{i}})| + |Tf(\lambda_{n_{i}})-Tf(\eta)| 
\le
$$
$$
\le \varepsilon+|f(\lambda_{n_{i}})-f(A\lambda_{n_{i}})|+\varepsilon = 2\varepsilon + |f(\lambda_{n_{i}})-f(\lambda_{n_{i}})- f(\frac{1}{n_{i}}[A^{n_{i}+1}\mu-A\mu])| \le
$$
$$
\le 2\varepsilon+\frac{1}{n_{i}}|f(A^{n_{i}+1}\mu - A\mu)| \le 2\varepsilon + \frac{2\|f\|}{n_{i}}.
$$
         
Поскольку $n_{i} \to \infty$ при $i \to \infty$, то  $|f(\eta)-f(A\eta)| \le 2\varepsilon$.  Так как 
$\varepsilon$  произвольно,  то $|f(\eta)-f(A\eta)|=0$.

Итак, для каждой $f \in B(X, \Sigma)$ выполняется равенство $f(\eta)=f(A\eta)$.  Множество $B(X, \Sigma) )$ тотально на   
$ba(X, \Sigma)$. Следовательно, $\eta=A\eta$.  При этом $\eta \in ba(X, \Sigma)$. 

Покажем, что $\eta \in S_{ba}$,
т.е. мера $\eta$ нормирована и положительна. 

Рассмотрим $\tau_{B}$-окрестность точки $\eta$
вида $V(\eta,f,\varepsilon)$, где $\varepsilon>0$  произвольно, $f \in B(X, \Sigma)$ и $f(x) \equiv 1$.
Тогда существует $n_{i}$, такое что $\lambda_{n_{i}} \in V(\eta,f,\varepsilon)$, т.е.
$$
|(f,\eta)-(f, \lambda_{n_{i}})|=|\eta(X)-\lambda_{n_{i}}(X)|<\varepsilon.
$$  
Поскольку $\lambda_{n}(X)=1$ \; для всех $n \in N$, то $|\eta(X)- 1|<\varepsilon$ \; при любом $\varepsilon>0$, откуда $\eta(X)=1$. 

 Предположим, что существует $E \in \Sigma$ и $r>0$, такие, что $\eta(E)=-r<0$. \\

Возьмем в качестве $f \in B(X, \Sigma)$ характеристическую функцию $f=\chi_{E}$ множества $E$ и число $\varepsilon=\frac{r}{2}$. Тогда $\langle f,\eta \rangle=\eta(E)=-r, \langle f,\lambda_{n}\rangle=\lambda_{n}(E)$ и

$$
|\langle f,\lambda_{n}\rangle-\langle f,\eta \rangle| = \lambda_{n}(E)+r \ge r>\varepsilon
$$
при всех $n \in N$.

Следовательно, $\lambda_{n} \notin V(\eta, f, \varepsilon)$ при $n \in N$, т.е. мера $\eta$ не является $\tau_{B}$-предельной для последовательности  $\{\lambda_{n}\}$. Полученное противоречие доказывает, что $\eta(E) \ge 0$  для всех $E\in \Sigma$.

Суммируя все предыдущие заключения, мы получаем, что $\mathfrak{M}{\{\lambda_{n}\}} \subset \Delta_{ba}$,  $\mathfrak{M}{\{\lambda_{n}\}} \ne \O$ и  $\mathfrak{M}{\{\lambda_{n}\}}$ -- $\tau_{B}$-компактно.

Теорема доказана. \medskip

Нам понадобится также ещё одна теорема из ~\cite{Zhd01}.

Обращаем внимание на то, что в этой теореме средние берутся для ЦМ, имеющих при каждом $n \in N$ различные начальные меры $\mu_{n}$. Это делается не ради обобщения, а для конкретного дальнейшего использования.

\begin{theorem}\label{t4.3}
 (\cite[теорема 7.3, стр. 78]{Zhd01}). Пусть на произвольном $(X,\Sigma)$ задана произвольная ЦМ, меры $\mu_{n} \in ba(X,\Sigma)$, $\mu_{n} \in S_{ba}$ и
$$
\lambda_{n}=\lambda_{n}^{\mu_{n}} = \frac{1}{n} \sum_{k=1}^{n}A^{k}\mu_{n}, \;\; n=1,2,\ldots .
$$
Тогда каждая $\tau_{B}$-предельная мера последовательности $\{\lambda_{n}\}$ является инвариантной для оператора $A$,  т.е.  $\mathfrak{M}{\{\lambda_{n}\}} \subset \Delta_{ba}$,  множество таких мер непусто, $\mathfrak{M}{\{\lambda_{n}\}} \ne \O$ и  $\tau_{B}$-компактно. 
\end{theorem}



\begin{remark}\label{r4.1} 
Если $\eta \in \mathfrak{M}{\{\lambda_{n}\}}$, то для всех $E \in \Sigma$, для которых существует $\lim \lambda_{n}(E)$, очевидно, что $\eta(E)=\lim \lambda_{n}(E)$.\medskip
\end{remark}

  Теперь мы готовы перейти к доказательству одной из основных наших теорем о том, что классическое условие $(D)$ эквивалентно нашему условию $(*)$.
Оказалось, что это можно сделать по той же общей схеме, но с существенными новыми деталями, что и доказательство теоремы 4.1 (\cite[теорема 12.2]{Zhd02}), приведенной выше.

Прежде чем привести эту новую теорему, мы докажем одну техническую лемму, которая нам понадобится при доказательстве второй части (необходимости) теоремы 4.4.

 \begin{lemma}\label{l4.1} 
 Для любой ЦМ, для любого множества $G \in \Sigma$ выполняется равенство
$$
\sup_{x\in X}p^{m}(x, G)= \sup_{\eta \in S_{ba}}A^{m}\eta(G).
$$
\end{lemma}

Подчеркнем, что супремум справа берется по всем конечно-аддитивным вероятностным мерам из $S_{ba}$, т. е. по всем счётно-аддитивным и по всем чисто конечно-аддитивным мерам $\eta \in S_{ba}$.

{\bf Доказательство.} 
Мы уже отмечали выше, что переходную функцию при любых $m \in N, x \in X,$ и $G \in \Sigma$ можно представить в следующем виде
$$
p^{m}(x, G)=A^{m}\delta_{x}(G),
$$
где $\delta_{x}$ - мера Дирака в точке $x \in X$.

Отсюда следует, что
$$
\sup_{x\in X}p^{m}(x, G)= \sup_{x\in X}A^{m}\delta_{x}(G).
$$

Тогда очевидно, что при увеличении класса мер, по которому ищется супремум, этот супремум может только возрасти:
$$
\sup_{x \in X}A^{m}\delta_{x}(G) \le \sup_{\eta \in S_{ba}}A^{m}\eta (G).
$$

Покажем, что на самом деле здесь имеет место равенство.
Предположим, что неравенство, записанное выше, строгое. Введём числа $\alpha, \beta > 0$ такие, что $\alpha < \beta$ и
$$
\sup_{x\in X}p^{m}(x, G)=\sup_{x \in X}A^{m}\delta_{x}(G) =\alpha<\beta = \sup_{\eta \in S_{ba}}A^{m}\eta (G).
$$

Тогда существует такая мера $\eta_{0} \in S_{ba}$, что 
$$
\alpha< \alpha+\frac{\beta-\alpha}{2} \le A^{m}\eta_{0} (G) \le \beta.
$$

Но, при этом
$$
A^{m}\eta_{0} (G) =\int_{X}p^{m}(x, G)\eta_{0}(dx) \le \int_{X}\alpha \cdot \eta_{0}(dx) =\alpha.
$$

Полученное неравенство противоречит предыдущему: $ A^{m}\eta_{0} (G) \ge \alpha+\frac{\beta-\alpha}{2}$. Следовательно, верно равенство в формулировке Леммы 4.1. Лемма доказана. \\
 
 \begin{theorem}\label{t4.4} 
 Для произвольной ЦМ условие Деблина $(D)$ эквивалентно условию (*):
$$ (*) \;\;\;\; \qquad \Delta_{ba}\subset{ca(X,\Sigma)}, $$
которое означает, что все инвариантные конечно-аддитивные меры исходной ЦМ являются счетно-аддитивными, или, другими словами, исходная ЦМ не имеет инвариантных чисто конечно-аддитивных мер.
\end{theorem}

{\bf Доказательство.}
  Вначале докажем, что из выполнения условия $(*)$ следует выполнение условия $(D)$ (достаточность).

Пусть выполнено условие $(*)$. Тогда по теореме 8.2 ~\cite{Zhd01} $dim\Delta_{ba}=dim\Delta_{ca}=n<\infty$. Пусть $H_{ba}=H_{ca}=\{\mu_{1}, ..., \mu_{n} \}$ сингулярный базис пространства инвариантных мер $M_{ba}=M_{ca}$, существующий по теореме 6.3 ~\cite{Zhd01}.

Введём счётно-аддитивную меру $\varphi=\mu_{1}+ ...+\mu_{n}$. Мы хотим доказать, что для этой меры $\varphi$ и при некоторых числах $m\ge 1 (m\in N)$ и $\varepsilon >0$ выполнено условие Деблина $(D)$.

Предположим, что для данной меры $\varphi$ не выполнено условие $(D)$. Тогда для каждого  $m\ge 1 (m\in N)$ и для каждого $0<\varepsilon<1$, существует множество  $E_{m, \varepsilon} \in \Sigma$ и существует точка $x_{m, \varepsilon} \in X$, такие, что
 $$
 \varphi (E_{m, \varepsilon}) < \varepsilon \;\;  \text{ и } \;\;  p^{m}(x_{m, \varepsilon}; E_{m, \varepsilon}) > 1-\varepsilon.
 $$
 
 Для каждого $m=1, 2, ...$ и для некоторго фиксированного $\delta \in (0, 1)$ положим $\varepsilon =\varepsilon(m) =\frac{\delta}{2^{m}}$, $E_{\delta} = \cup_{m=1}^{\infty}E_{m, \epsilon(m)}$.
  
  Тогда, так как мера  $\varphi$ счётно-аддитивна, а значит и счётно-полуаддитивна, то
 $$
 \varphi (E_{\delta}) =\varphi (\bigcup_{m=1}^{\infty}E_{m, \epsilon(m)}) \le \sum_{m=1}^{\infty} \varphi (E_{m, \varepsilon(m)}) \le \sum_{m=1}^{\infty}\varepsilon(m)=\delta \sum_{m=1}^{\infty}\frac{1}{2^{m}} = \delta,
 $$
 $\text{ т. е. } \varphi(E_{\delta})<\delta.$
 
При этом, для всех $m=1, 2, ...$ выполняется:
$$
 p^{m}(x_{m, \varepsilon(m)}; E_{\delta}) \ge p^{m}(x_{m, \varepsilon(m)}; E_{m, \varepsilon(m)}) > 1-\varepsilon(m) = 1 - \frac{\delta}{2^{m}}.
$$

Теперь построим последовательность мер Дирака $\eta_{m} \in S_{ca}$, сосредоточенных в точках $x_{m, \varepsilon(m)}$, т. е. $\eta_{m} = \delta_{x_{m, \varepsilon(m)}}$. Тогда для степеней марковских операторов $A^{s}$ при $s=1, 2, ..., m$ будет выполняться:
$$
A^{s} \eta_{m}(E) =  p^{s}(x_{m, \varepsilon(m)}; E) 
$$
для произвольных множеств $E \in \Sigma$.

Соответственно для средних по Чезаро при $E=E_{\delta}$ и $m=1, 2, ...$ будет верно следующее:
$$
 \lambda_{m}^{\eta_{m}}(E_{\delta}) = \frac{1}{m}\sum_{s=1}^{m}A^{s}\eta_{m}(E_{\delta})=
  \frac{1}{m}\sum_{s=1}^{m}p^{s}(x_{m, \varepsilon(m)}; E_{\delta}) >
$$    
$$
  >  \frac{1}{m}\sum_{s=1}^{m}(1 - \frac{\delta}{2^{s}}) = 
        \frac{1}{m}(m - \sum_{s=1}^{m}\frac{\delta}{2^{s}})=
          1- \frac{\delta}{m}(1 - \frac{\delta}{2^{m}}) > 1- \frac{\delta}{m},
$$
$
\text{ т. е. }  \lambda_{m}^{\eta_{m}}(E_{\delta}) > 1- \frac{\delta}{m}, m=1, 2, ... .
$

Теперь мы обратимся к $\tau_{B}$-слабой топологии в пространстве $ba(X, \Sigma)$, порождённой предсопряжённым пространством $B(X, \Sigma)$. По теореме 4.3 (\cite[теорема 7.3]{Zhd01} множество $\mathfrak{M} \{ \lambda_{m}^{\eta_{m}} \}$ всех $\tau_{B}$-предельных точек (мер) последовательности $\lambda_{m}^{\eta_{m}}$ непусто и содержится в множестве всех инвариантных для ЦМ конечно-аддитивных мер,
$\mathfrak{M} \{ \lambda_{m}^{\eta_{m}} \} \subset \Delta_{ba} = \Delta_{ca}$, т. е. все
$\tau_{B}$-предельные точки (меры) последовательности $\{ \lambda_{m}^{\eta_{m}} \}$ инвариантны и, по условию $(*)$, счётно\--ад\-ди\-тив\-ны.

Пусть мера $\mu \in \mathfrak{M} \{ \lambda_{m}^{\eta_{m}} \}$ произвольна. Тогда, если для некоторого  $E \in \Sigma$ существует предел числовой последовательности $ \lambda_{m}^{\eta_{m}}(E)$ при
$n \to \infty$, то, очевидно, $\mu(E)=\lim_{n \to \infty} \lambda_{m}^{\eta_{m}}(E).$
Возьмём множество $E=E_{\delta}$. Тогда из полученной выше оценки $\lambda_{m}^{\eta_{m}}(E_{\delta}) > 1- \frac{\delta}{m}, m=1, 2, ... $, следует, что $\lambda_{m}^{\eta_{m}}(E_{\delta}) \to 1$ при $m \to \infty$. Следовательно, для $\tau_{B}$-предельной (инвариантной) меры $\mu$ выполняется $\mu(E_\delta)=1$. При этом, как мы получили выше, $\varphi(E_{\delta}) < \delta < 1$. 

Поскольку мера $\mu$ инвариантна, то она разложима по сингулярному базису
$H_{ba}=\{\mu_{1}, ... , \mu_{n}\}, \mu=\sum_{i=1}^{n} \alpha_{i}\mu_{i}$,  где $0 \le \alpha_{i} \le 1, i=1, 2, ..., n$. 
Отсюда, $\mu=\sum_{i=1}^{n} \alpha_{i}\mu_{i} \le \sum_{i=1}^{n}\mu_{i}=\varphi$.
Следовательно, $\mu(E_{\delta}) \le \varphi(E_{\delta}) < \delta < 1.$
Но выше мы получили, что $\mu(E_{\delta}) =1$.
Возникшее противоречие доказывает, что условие 
$(*)$ влечёт выполнение условия $(D)$ для меры $\varphi=\mu_{1}+\mu_{2}+ ... + \mu_{n}$.
Первая часть теоремы доказана.\medskip

Теперь докажем, что из выполнения условия $(D)$ следует выполнение условия $(*)$ (необходимость).

Пусть выполнено условие $(D)$ для некоторых $\varphi \in ca(X,\Sigma), \varphi \ge 0,$ $\varepsilon >0$ и $m \ge 1 (m \in N)$.

Предположим, что условие $(*)$ не выполнено. Тогда ЦМ имеет инвариантную чисто конечно-аддитивную меру $\lambda \in S_{ba}, \lambda = A\lambda$.

Любая чисто конечно-аддитивная мера дизъюнктна с любой счётно-ад\-ди\-тив\-ной мерой (\cite[теорема 1.16]{YoHew1}),  откуда $\lambda \wedge \varphi =0$, т. е. меры  $\lambda$ и $\varphi$ дизъюнктны. Тогда, по другой теореме  (\cite[теорема 1.21]{YoHew1}), для любого числа $\varepsilon > 0$ (а значит и для нашего $\varepsilon$ из ``тройки'' 
($\varphi, \varepsilon, m$)) существует множество  $G \in \Sigma$ такое, что
$\varphi (G)<\varepsilon$ и $\lambda (X \setminus G)< \frac{\varepsilon}{2}$ (нам дальше нужно именно $\frac{\varepsilon}{2}$ во втором неравенстве, а не  $\varepsilon$, что позволяется теоремой 1.21 из \cite{YoHew1}). Отсюда  $\lambda (G) \ge 1- \frac{\varepsilon}{2}$.

Теперь воспользуемся леммой 4.1. Сделаем следующие преобразования и оценки:
$$
\sup_{x\in X}p^{m}(x, G)= \sup_{\eta \in S_{ba}}A^{m}\eta(G) \ge A^{m}\lambda(G) = \lambda(G)>1- \frac{\varepsilon}{2}.
$$
Следовательно, найдется хотя бы одна точка $x_{0} \in X$, такая, что $$
\sup_{x\in X}p^{m}(x, G) \ge p^{m}(x_{0}, G)>1- \frac{\varepsilon}{2} >1- \varepsilon.
$$
Но по условию $(D)$, которое предполагаем выполненным, должно быть $p^{m}(x, G)<1- \varepsilon$ для всех $x\in X$, в том числе и для $x=x_{0}$. Полученное противоречие доказывает, что условие $(*)$ также  выполнено. Теорема доказана.\medskip

Привлекательность условия (*) в том, что в нем, в отличие от условия Деблина $(D)$ и нашего условия $(\tilde D)$, нет никакой аналитики и динамики, а есть только «качественное» понятное утверждение.\medskip

Введём ещё одно условие:

$$
(\tilde*)   \;\;\;\; \Delta_{pfa} = \O.
$$

Выполнение условия $(\tilde*)$ означает, что у ЦМ нет инвариантных чисто конечно\- -ад\-ди\-тив\-ных мер, т. е. все её инвариантные конечно\- -ад\-ди\-тив\-ные меры (а такие всегда есть) являются счётно-ад\-ди\-тив\-ными.

 \begin{corollary}\label{cor4.1} 
  Для произвольной ЦМ условие $(*)$ эквивалентно условию $(\tilde*)$.
 \end{corollary}
 
\begin{corollary}\label{cor4.2} 
 Для произвольной ЦМ условие $(D)$ эквивалентно условию $(\tilde*)$.
 \end{corollary}
 
\begin{corollary}\label{cor4.3}
 Если для произвольной ЦМ выполнено условие $(*)$, то выполнено условие $(K)$, т. е. ЦМ квазикомпактна.
\end{corollary}
 
\begin{corollary}\label{cor4.4}
Если ЦМ не является квазикомпактной, то ЦМ имеет инвариантные чисто конечно-аддитивные меры, т. е. $\Delta_{pfa} \ne \O$. \medskip
 \end{corollary}
 
Итак, если для ЦМ выполнено условие $(*)$, то выполнено условие $(D)$ и выполнено условие квазикомпактности $(K)$. Следовательно, для такой ЦМ справедливы все эргодические теоремы, в том числе равномерные, хорошо и подробно изложенные в литературе в разных вариантах, в том числе с учётом циклических мер и с оценками скорости сходимости (см., например, книги Дж. Дуба \cite[глава V, \S 5]{Doob1}, М. Лоэва \cite[глава IX, \S 32]{Loeve1}, Ж. Невё \cite[глава V, пункт V.3]{New1}). 

Мы не собираемся повторять здесь эти объёмные известные теоремы, но один самый короткий вариант из них мы всё-таки воспроизведём, поскольку он понадобится нам и в следующем \S 5.

Для этого воспользуемся книгой ~\cite{DS1}, где изучение квазикомпактных ЦМ проводится на более близком нам операторном языке и без подробностей.   
В теореме VIII.8.3 (\cite{DS1}, стр. 754) и в следствии VIII.8.4. ~\cite{DS1} доказывается нужное нам утверждение (при условии, что $A^{n}/n$ сходится к нулю в слабой топологии, которое в нашем случае выполнено, так как $\|A^{n}\|=1$ при $n\in N$).
Из указанной теоремы получаем как простое следствие следующее утверждение.

\begin{theorem}\label{t4.5} 
(равномерная эргодическая теорема для условия $(*)$)
Пусть для произвольной ЦМ выполнено условие $(*)$, 
которое означает, что все инвариантные конечно-аддитивные меры исходной ЦМ являются счётно\--ад\-ди\-тив\-ными, или, другими словами, исходная ЦМ не имеет инвариантных чисто конечно-аддитивных мер.

Тогда пространство $M_{ca}$ инвариантных счётно-аддитивных мер ЦМ конечномерно, т. е. $dim \Delta_{ba}=dim \Delta_{ca}=n<\infty$,  ему  
соответствует конечномерный проектор  $E(1, A)$ с собственным значением $\lambda = 1$, и
 $$\frac{1}{n} \sum_{i=1}^{n} A^{i} \Rightarrow E(1, A), n \to \infty,$$
 т. е. последовательность средних по Чезаро от степеней оператора $A$ сходится к проектору $E(1, A)$ в равномерной топологии операторов (т. е. равномерно по всем начальным  счётно-аддитивным мерам).
 \end{theorem} 

\begin{remark}\label{r4.2} 
Доказательство теорем в ~\cite[глава VIII, \S 8]{DS1} проводится в рамках спектральной теории линейных операторов с использованием комплексного расширения исходного банахова пространства, однако в формулировках используемых нами теоремы 8.3 из ~\cite[глава VIII, \S 8]{DS1} и следствия 8.4 из \cite[глава VIII, \S 8]{DS1} комплексные числа уже не используются.
\end{remark}

\begin{remark}\label{r4.3} 
В ~\cite{DS1} термины ``компактный'' (т. е. вполне непрерывный) и ``квазикомпактный'' оператор не используются.
\end{remark}

\begin{remark}\label{r4.4} 
В начале \S 8 Главы VIII работы ~\cite{DS1} говорится, что результаты этого параграфа ``могут быть успешно применены к некоторым задачам в теории марковских процессов, хотя мы лишь бегло укажем эти применения''.
\end{remark}

Эти применения приводятся в конце этого же параграфа 8 Главы VIII работы ~\cite{DS1}, где вводится два марковских оператора:
$T: B(X, \Sigma) \to B(X, \Sigma)$ и двойственный (но не сопряжённый) к нему оператор $T ': ca(X, \Sigma) \to ca(X, \Sigma)$, интегрально определяемые через переходную функцию ЦМ.

Однако, продолжение оператора $T '$ на пространство $ba(X, \Sigma)$ там не производится, и конечно-аддитивные меры не используются. В этом и заключается главное отличие развиваемого нами операторного подхода в теории классических цепей Маркова, от операторного подхода, предлагаемого в книге ~\cite{DS1}.

Из доказанной теоремы 4.4 и из теоремы 4.1 (\cite[теорема 12.2]{Zhd02}) мы, как следствие, получаем следующее обещанное и психологически важное для нас утверждение.

\begin{theorem}\label{t4.6}
Для произвольной цепи Маркова условие $(D)$ эквивалентно условию $(\tilde{D})$.
\end{theorem} 

Эта теорема позволяет при дальнейшем цитировании некоторых утверждений из ~\cite{Zhd01} и ~\cite{Zhd02} заменять маркер $(\tilde{D})$ на маркер $(D)$, что мы и будем делать без упоминания об этом.

\section{Теоремы о размерности множества инвариантных мер и квазикомпактности цепи Маркова}\label{s5}

В работе автора ~\cite{Zhd01} была доказана (без использования условия квазикомпактности) теорема 8.1 о том, что если $dim\Delta_{ca}=\infty$, то множество $\Delta_{pfa}$ непусто и $dim\Delta_{pfa}=\infty$.  Напомним, что $\Delta_{pfa}$ -- это множество всех инвариантных вероятностных чисто конечно-аддитивных мер. Из этой теоремы сразу же  вытекает следующее утверждение (см. ~\cite[теорема 8.2]{Zhd01}), которое мы будем использовать.

\begin{theorem}\label{t5.1}
(см. \cite[теорема 8.2]{ Zhd01}). 
Для произвольной цепи Маркова, если выполнено условие (*), т.е. если $\Delta_{ba}\subset{ca(X,\Sigma)}$, то выполнено условие
$$  (**) \;\;\; dim\Delta_{ba}=n<\infty,$$
$$\;\;\;\;  \text {т. е. } (*) \Rightarrow (**). $$
\end{theorem}

{\bf Следствие 5.1}. 
 Если для произвольной цепи Маркова выполнено условие $(*)$, то $\Delta_{ba}=\Delta_{ca} \ne \O$.

{\bf Следствие 5.2}. 
Если для произвольной цепи Маркова выполнено условие $(D)$ или $(\tilde D)$, то выполнено условие $(**)$.







Возникает интуитивное предположение, что утверждение теоремы 5.1 (~\cite[теорема 8.2]{Zhd01}) можно обратить. В той же работе ~\cite[теорема 8.3]{Zhd01} было доказано такое обращение, но только для случая размерности $n=1$.
 Эта теорема имеет для нас важное значение, и мы приведём её ниже в усиленном виде (по сравнению с теоремой 8.3 из ~\cite{Zhd01}). 
 
 Пусть $dim \Delta_{ba}=1$, $\Delta_{ba}=\{ \mu \}$, т. е. оператор $A=A^{1}$ имеет единственную инвариантную конечно-аддитивную меру $\mu \in S_{ba}$. Тогда возможны два взаимоисключающих случая.
 
 Первый - это когда $\mu$ является взвешенной суммой циклических мер одного цикла $\{ \mu_{1}, \mu_{2}, ..., \mu_{n} \}$, $\mu = \frac{1}{n} \sum_{i=1}^{n}\mu_{i}$, $A\mu_{i}=\mu_{i+1}, A\mu_{n}=\mu_{1}, n \ge 2$, где все $\mu_{i}$ - это инвариантные меры оператора $A^{n}$, $\mu_{i} = A^{n}\mu_{i}, i = 1, 2, ..., n$. В этом случае $\mu = A^{1}\mu$. Назовём такую инвариантную меру $\mu$ {\it составной}. Заметим, что период этого цикла $n$ может быть любым ($n \ge 2$).

Второй случай - это когда $\mu$ не является взвешенной суммой циклических мер, т. е. у ЦМ нет циклов ни конечно-аддитивных, ни счётно-аддитивных. Назовём такую инвариантную меру $\mu$ {\it простой}.
В этом случае соответствующую цепь Маркова назовём, по аналогии со счётно-аддитивным случаем, {\it апериодической}.

Теперь докажем усиленный вариант теоремы 8.3 из ~\cite{Zhd01}.

\begin{theorem}\label{t5.2}
Пусть на  некотором $(X,\Sigma)$ задана произвольная ЦМ. Если $dim \Delta_{ba}=1$, т. е. если ЦМ  имеет в 
$S_{ba}$ единственную инвариантную меру $\mu$, $\Delta_{ba}=\{ \mu \}$, то выполнено условие $(*)$:  $\Delta_{ba} \subset ca(X,\Sigma)$, т.е.  инвариантная 
мера $\mu$  счетно-аддитивна. При этом для любой начальной меры $\eta \in S_{ba}$ имеет место сходимость $ \lambda_{n}^{\eta} \to \mu$  в $\tau_{B}$-топологии, т. е. $\lim_{n \to \infty} \lambda_{n}^{\eta}(E) = \mu (E)$, для всех $E \in \Sigma$.

Если же инвариантная мера $\mu$ простая, то не только средние по Чезаро $\{\lambda_{n}^{\eta}\}$ сходятся к $\mu$, но и сама марковская последовательность мер $\mu_{n} = A\mu_{n-1} = A^{n}\eta$ сходится к $\mu$ в $\tau_{B}$-топологии при любой начальной мере $\eta \in S_{ba}$, т. е.  $\mu_{n}(E) \to \mu(E)$  для всех $E \in \Sigma$ при $n \to \infty$.
\end{theorem} 

{\bf Доказательство.}
Пусть $\eta \in S_{ba}$ и $\mathfrak M$ -- множество всех $\tau_{B}$-предельных мер последовательности 
$\{\lambda_{n}^{\eta}\}$.  По теореме 4.2 имеем $\mathfrak{M} \ne \O$  и  $\mathfrak{M} \subset \Delta_{ba}$. Отсюда и из условий теоремы следует, что $\mathfrak{M}=\Delta_{ba}=\{\mu\}$, т.е. $\{\lambda_{n}^{\eta}\}$ имеет единственную 
$\tau_{B}$-предельную точку $\mu$. 

Покажем, что $\lambda_{n}^{\eta} \to \mu$ в $\tau_{B}$-топологии. 
Пусть это не так, т.е. существуют $E \in \Sigma$, $\delta>0$ и $\{n_{i}\}$ такие, что 
$|\lambda_{n_{i}}^{\eta}(E)-\mu(E) | \ge \delta$, $i=1,2,\ldots .$ Тогда подпоследовательность 
$\{\lambda_{n_{i}}^{\eta}\}$  имеет $\tau_{B}$-предельную точку $\xi \ne \mu$. Но все $\tau_{B}$-предельные точки для любой подпоследовательности  $\{\lambda_{n_{i}}^{\eta}\}$  
являются $\tau_{B}$-предельными и для самой последовательности $\{\lambda_{n}^{\eta}\}$,  т.е. 
$\xi \in \mathfrak{M}=\{\mu\}$.  Полученное противоречие показывает, что $\lambda_{n}^{\eta} \to \mu$   в 
$\tau_{B}$-топологии для любой $\eta \in S_{ba}$, т. е. $\lambda_{n}^{\eta}(E) \to \mu(E)$ для всех $E \in \Sigma$ при $n \to \infty$.  

Пусть теперь начальная мера $\eta$ счетно-аддитивна, т.е.  $\eta \in S_{ca}$. Тогда все 
$\lambda_{n}^{\eta}$, \; $n=1,2,\ldots$ также счетно-аддитивны, и, как мы уже показали, $\lambda_{n}^{\eta}(E) \to \mu(E)$ для любого $E \in Z$. Тогда, по теореме  Никодима  
(\cite[глава III, \S 7, cледствие 4]{DS1}), отсюда вытекает счетная аддитивность меры $\mu$. 

В книге Дж. Дуба (\cite[глава V, \S 5, случай $e)$, стр. 196]{Doob1}) утверждается, что, если ЦМ не имеет циклов, то предел (в некотором смысле) средних по Чезаро от степеней (свёрток) переходных вероятностей ЦМ совпадает с пределом от самих этих степеней. Если перейти к нашим операторным терминам, то это означает, что вместо пределов последовательностей мер $\{\lambda_{n}^{\eta}\}$
 можно использовать равные им пределы самой марковской последовательности мер $\mu_{n}=A\mu_{n-1}=A^{n}\eta, n \in N$.
 При этом, очевидно, тип сходимости (сильная, слабая) здесь не важен.
 
 Пусть наша единственная инвариантная мера $\mu=A\mu$ простая, что по определению означает, что у ЦМ нет циклов. Тогда из приведённого утверждения Дуба следует, что $\mu_{n}(E) \to \mu(E)$ для всех $E \in \Sigma$ при $n \to \infty$, т. е. 
 $\mu_{n} \to \mu$   в $\tau_{B}$-топологии для любой начальной меры $\eta \in S_{ba}$.
Теорема доказана. \medskip



\begin{remark}\label{r5.1} 
Подчеркнем, что в теореме 5.2 начальная мера $\eta \in S_{ba}$ может быть не только счётно-аддитивной, но и чисто конечно-аддитивной.
\end{remark}

\begin{remark}\label{r5.2} 
Далее, в теореме 5.4 мы покажем, что сходимость последовательностей мер $\lambda_{n}^{\eta}$ и $\mu_{n}$ в теореме 5.2 не только слабая, но и сильная метрическая. 
\end{remark}

{\bf {Следствие 5.3}.} 
Если $dim \Delta_{ba}=1$, то у ЦМ нет инвариантных чисто конечно-аддитивных мер, т. е.  выполнено условие $(\tilde{*})$.

Пусть $dim\Delta_{ba}=1, \Delta_{ba}=\{\mu\}$. Тогда, по теореме 5.2, 
$\mu \in ca(X, \Sigma)$, т. е. мера $\mu$ счётно-аддитивна, и  $\Delta_{ba}=\Delta_{ca} \subset ca(X, \Sigma)$. Следовательно, выполнено условие $(*)$. По теореме 4.4 условие $(*)$ эквивалентно условию Деблина $(D)$.


Отсюда вытекает следующая теорема.

\begin{theorem}\label{t5.3}
Пусть для произвольной цепи Маркова $dim\Delta_{ba}=1$. Тогда выполнены условия $(D)$, $\tilde D$ и $(K)$, т. е. ЦМ квазикомпактна, причем $ \Delta_{ba}=  \Delta_{ca}=\{\mu\}\subset ca(X, \Sigma)$, т. е. выполнено условие $(*)$.
\end{theorem} 
 
\begin{theorem}\label{t5.4} 
(равномерная эргодическая теорема для $n=1$).
Пусть для некоторой цепи Маркова $dim\Delta_{ba}=1, \Delta_{ba}=\{\mu\}$, $\mu=A\mu$.
Тогда расстояние 
$$\rho({\lambda_{n}^{\eta}}, \mu)=\|{\lambda_{n}^{\eta}} - \mu \| \to 0 \text{ при } n\to\infty$$ 
равномерно  по начальным мерам $\eta \in S_{ca}$, т. е. последовательность мер ${\lambda_{n}^{\eta}}$ равномерно сходится к инвариантной мере 
$\mu$ в сильной метрической топологии.

Если при этом мера $\mu$ является простой, т. е. у ЦМ нет циклов, то  
$$\rho(\mu_{n}, \mu)=\|\mu_{n} - \mu \| \to 0 \text{ при } n \to\infty$$ 
равномерно  по начальным  мерам $\eta \in S_{ca}$, где $\mu_{n}=A^{n}\eta$, $n\in N$.
\end{theorem}

{\bf Доказательство.}
По теореме 5.2 мера $\mu$ счётно-аддитивна. По теореме 5.3 данная ЦМ квазикомпактна. Следовательно, для неё верны все замечания, предшествующие теореме 4.5, и верно само утверждение  теоремы 4.5. При этом проектор  $E(1, A)$ одномерен. Поэтому первое утверждение данной теоремы является частным случаем общего утверждения теоремы 4.5.

Второе утверждение для простой инвариантной меры $\mu$ получаем из тех же утверждений Дж. Дуба ~\cite{Doob1}, что и в доказательстве  теоремы 5.2. При этом мы усиляем утверждение теоремы 5.2 о слабой сходимости соответствующих последовательностей мер, гарантируя их сильную сходимость равномерно по всем начальным счётно-аддитивным мерам.
Теорема доказана.

 В настоящей работе мы приводим и доказываем обращение теоремы 5.1, т.е. обобщаем и усиляем теорему 5.2 уже для произвольной размерности $n\in{N}$, но при дополнительных условиях $(\alpha)$ и $(\beta)$:

$ 
(\alpha)
\begin{cases}
\text{{\it Пусть $\mu$ инвариантная конечно-аддитивная мера, $\mu\in\Delta_{ba}$,}}\\
\text{{\it и пусть дано множество $K_{\mu}\in\Sigma$, такое, что  $\mu(K_{\mu})=1$.}}\\
\text{{\it Тогда существует множество $K\subset{K_{\mu}}$, такое, что $\mu(K)=1$ и }}\\
\text{{\it оно стохастически замкнуто, т.е. $p(x,K)=1$ для любого $x\in{K}$.}}
\end{cases}
$ \medskip

$ 
(\beta)
\begin{cases}
\text{{\it Пусть пространство $M_{ba}$  всех инвариантных конечно-ад\-ди\-тив\-ных }}\\
\text{{\it мер ЦМ конечномерно и $H_{ba} = \{ \mu_{1}, \mu_{2}, ..., \mu_{m} \}$ его дизъюнктный}}\\
\text{{\it базис, существующий по теореме 3.3. Тогда инвариантные меры}}\\
\text{{\it из $H_{ba}$ попарно сингулярны.}}
\end{cases}
$ \medskip

Очевидно, что при  $K_{\mu}=X$ условие $(\alpha)$ тривиально выполняется при $K=K_{\mu}=X$.

Если предположить, что мера $\mu\in\Delta_{ba}$ счётно-аддитивна, то условие $(\alpha)$ также будет выполнено (\cite[теорема 6.2]{Zhd01}). \medskip

Если предположить, что все меры из $H_{ba}$ счётно-аддитивны, то условие $(\beta)$ выполнено. Пока что по теореме 3.3 все меры из $M_{ba}$ всего лишь попарно дизъюнктны. 

\begin{theorem}\label{t5.5}
Рассматриваем произвольную цепь Маркова на общем фазовом пространстве $(X,\Sigma)$. Пусть пространство её инвариантных конечно-ад\-ди\-тив\-ных мер конечномерно (выполнено условие $(**)$), т. е. $dimM_{ba}=dim\Delta_{ba}=n<\infty$, и выполнены условия $(\alpha)$ и $(\beta)$. 

Тогда  $\Delta_{ba}\subset{ca(X,\Sigma)}$, т. е. все инвариантные конечно-аддитив\-ные меры цепи Маркова являются счетно-аддитивными, выполнены условия $(*)$, эквивалентные ему условия $(D)$ и $(\tilde D)$, и ЦМ является квазикомпактной (т. е. выполнено условие $(K)$).
\end{theorem} 

{\bf Доказательство.}
Пусть $H=\{ \mu_{1}, ... , \mu_{n} \}$ - дизъюнктный базис в $M_{ba}$, существующий по теореме 3.3.
Из попарной сингулярности (условие $(\beta)$) конечного семейства мер $H=\{ \mu_{1}, ... , \mu_{n} \}$, следует, что существуют попарно непересекающиеся множества $K_1,K_2,$ $\dots,K_n\in\Sigma$ такие, что $\mu_i(K_i)=1$ для  $i=1,2,\dots,n$, и $K_i\cap{K_j}=\O$ при $i\ne{j}$.

По условию $(\alpha)$, для каждой меры $\mu_i$, $i=1,\dots,n$, существует множество $K^i\subset{K_i}$ такое, что $\mu_i(K^i)=1$ и $p(x,K^i)=1$ для всех $x\in{K^i}$, т. е. все $K^i$ стохастически замкнуты.

Образуем подпространства $(X_i,\Sigma_i)$ пространства $(X,\Sigma)$ по правилу $X_i=K^i\cap{X}$, $\Sigma_i=K^i\cap\Sigma$ для $i=1,2,\dots,n$. Сузим все меры $\mu_i$ на подпространства  $(X_i,\Sigma_i)$.

Сузим переходную функцию $p(x,E)$ исходной ЦМ на каждое подпространство $(X_i,\Sigma_i)$ по тождественному правилу: $\forall{x}\in{K^i}$, $\forall{E}\in\Sigma_i$, $p_i(x,E)=p(x,E)$. Легко проверить, что все меры $\mu_i$ будут единственными инвариантными конечно-аддитивными мерами для частных цепей Маркова, порождаемых переходными функциями $p_i(x,E)$ на $(X_i,\Sigma_i)$ при $i=1,2,\dots,n$.

Можно считать, что на каждом измеримом пространстве $(X_i,\Sigma_i)$ задана своя цепь Маркова с единственной инвариантной конечно-ад\-ди\-тив\-ной мерой $\mu_i$. Следовательно, по теореме 5.4, все инвариантные меры $\mu_i$ счетно-аддитивны на измеримых пространствах $(X_i,\Sigma_i)$ соответственно.

Теперь продолжим эти меры $\mu_i$ до мер $\widetilde{\mu_i}$ на всё исходное измеримое пространство $(X,\Sigma)$ по правилу $\widetilde{\mu_i}(E)=\mu_i(E\cap{K^i})$ для каждого $E\in\Sigma$, т.е. "нулями". Очевидно, счетно-аддитивность мер $\widetilde{\mu_i}$ при этом сохраняется и они будут инвариантны уже на всем $(X,\Sigma)$. 

По построению, $\widetilde{\mu_i}=\mu_{i}$ при всех $i=1,2,\dots,n$. Следовательно, все меры $\mu_i \in H$ счётно-аддитивны, и $\Delta_{ba}=\Delta_{ca} \subset 
ca(X,\Sigma)$, т. е. выполнено условие $(*)$  из теоремы 4.1 и выполнены эквивалентные ему условия $(D)$ и $(\tilde D)$, и ЦМ является квазикомпактной, т. е. выполнено условие $(K)$. Теорема доказана.\medskip

{\bf {Следствие 5.4}} 
(условное обращение теоремы 5.1).
 Для произвольной ЦМ при выполнении условий $(\alpha)$ и $(\beta)$  из условия $(**)$ следует выполнение условия $(*)$.

Случай $dimM_{ba} = 1$ был рассмотрен нами в теоремах 5.2, 5.3 и 5.4. Для таких ЦМ условие $(\alpha)$ выполнено автоматически при $K_{\mu}=K=X$, а условие $(\beta)$ теряет свой смысл. Поэтому в этих теоремах данные условия и не использовались.

\begin{theorem}\label{t5.6} 
(равномерная эргодическая теорема для $n \ge 2$). 
Пусть для произвольной цепи Маркова выполнены условия теоремы 5.5, т. е. выполнены условия $(**), dim M_{ba}=n, 2 \le n \le \infty,$ а также условия $(\alpha)$ и $(\beta)$. 

Тогда для цепи Маркова $M_{ba}=M_{ca}$ и верна равномерная эргодическая теорема в форме теоремы 4.5: последовательность средних по Чезаро от степеней марковского оператора $\frac{1}{m} \sum_{i=1}^{m}A^{i}$ сходится к конечномерному проектору на пространство инвариантных счётно-аддитивных мер $E(1, A)$ в равномерной топологии операторов.
\end{theorem}

{\bf Доказательство.}
По теореме 5.5 все инвариантные меры данной ЦМ являются счётно-аддитивными и ЦМ является квазикомпактной.
Следовательно, верно утверждение теоремы 4.5, что и записано в утверждении настоящей теоремы. Теорема доказана.\medskip

Приведя в статье эргодические теоремы 4.5, 5.4 и 5.6, мы преследовали и одну дополнительную цель.
Показать, что несмотря на использование в настоящей работе на каждом шагу ``экзотических'' конечно-аддитивных мер, финальные утверждения в эргодических теоремах уже их не содержат, т.е. они имеют классический характер. Вот только дорога к этим эргодическим теоремам была у нас другая -- через характеристики цепей Маркова, выраженные в терминах конечно-аддитивных мер.

В препринте ~\cite{Zhd06} мы показали, насколько этот аппарат расширяет набор наших инструментов при изучении внешне весьма простых примеров классических цепей Маркова даже на отрезке [0, 1].\medskip

{\bf Справка}

В работе Ш. Горовица ~\cite{Ho1} (1972 год) решаются некоторые близкие к нашей тематике задачи, но в другой постановке. Там рассматриваются цепи Маркова, заданные на измеримом пространстве $(X, \Sigma, m)$ с уже заранее заданной фиксированной ограниченной счётно-аддитивной мерой $m$. Переходная вероятность цепи Маркова $p(x,\cdot)$ предполагается при каждом $x \in X$ абсолютно непрерывной относительно меры $m$. Порожденные функцией $p(x, E)$ марковские операторы $P$ действуют (слева и справа) в пространствах $L_{1}(X, \Sigma, m)$ и в  $L_{\infty}(X, \Sigma, m)$, соответственно.

В теореме 4.1 ~\cite{Ho1} доказывается, что, если цепь Маркова эргодическая и консервативная (определения даны в начале статьи ~\cite{Ho1} -- не будем их здесь повторять), то квазикомпактность марковского оператора $P$ на пространстве $L_{\infty}(m)$ (условие $\{h\}$) эквивалентна отсутствию у него инвариантной чисто конечно-аддитивной меры (условие $\{a\}$) (pure charge). При этом оператор $P$ имеет единственную инвариантную счётно-аддитивную меру $\mu=\mu P$ (условие $\{d\}$), и цепь Маркова является харрисовской (условие $\{g\}$).

В книге Д. Ревюза (см. \cite[глава 6, §3, теоремы 3.5 и 3.7, стр. 240--243 в русском издании]{Rev1}) приводится в преобразованном виде указанная выше теорема Горовица ~\cite{Ho1} с новыми добавлениями и условиями. В частности, к общим условиям теоремы 4.1 ~\cite{Ho1} добавляется общее требование, чтобы цепь Маркова была возвратной по Харрису (см. ~\cite[теорема 6.3.7]{Ho1}), а мера $m$  из априорного пространства $(X, \Sigma, m)$ заранее объявляется инвариантной (см. \cite[глава 3, §2, определение 2.6]{Rev1}).

При этих жестких предположениях в ~\cite{Rev1} доказывается, что условие – {\it цепь Маркова является квазикомпактной}, эквивалентно условию – {\it не существует инвариантной чисто конечно-аддитивной меры} для цепи Маркова (см. \cite[теоремы 3.5 и 3.7]{Rev1}). Утверждается, что при выполнении этих эквивалентных условий заранее заданная инвариантная счетно-аддитивная мера $m$ из условия Харриса оказывается ограниченной и {\itединственной} инвариантной конечно-аддитивной мерой для цепи Маркова (\cite[теорема 3.7]{Rev1}).

В наших рассмотрениях цепь Маркова не предполагается харрисовской, и она вообще не привязана к какой-либо заранее заданной и уже инвариантной мере $m$. Наши теоремы 4.4, 5.3, 5.4, 5.5, следствие 4.3,  а также теоремы 8.2, 8.3 из ~\cite{Zhd01} и теорема 12.2 из ~\cite{Zhd02} носят более общий характер и содержат более сильные утверждения, чем в работах ~\cite{Ho1} и ~\cite{Rev1}.
\medskip 

{\bf Приложения.}

Укажем на возможности применения наших основных теорем в некоторых физических задачах на примере статьи ~\cite{Cio2019}. Её авторы разрабатывают специальный математический аппарат, называемый ими ``Термодинамический формализм'', для изучения соответствующих физических процессов. В частности, строится специальная классическая цепь Маркова на некотором бесконечномерном топологическом пространстве, но в качестве инвариантных мер привлекаются и конечно-аддитивные меры (чисто конечно-аддитивные) -- это используемая нами конструкция. В статье ~\cite{Cio2019} доказывается, что инвариантная вероятностная конечно-аддитивная мера для данной специальной ЦМ существует, она единственна, и у неё нет чисто конечно-аддитивной компоненты. На основании этих фактов доказывается ``асимптотическая устойчивость'' соответствующих марковских операторов, т. е. сходимость марковской последовательности счётно-аддитивных мер к инвариантной счётно-аддитивной мере в некоторой специальной метрике.

Мы видим здесь точки соприкосновения разрабатываемой в ~\cite{Cio2019} конструкции с нашими результатами из статей ~\cite{Zhd01} и ~\cite{Zhd02} (и из настоящей статьи), поскольку основные теоремы из указанных статей доказаны для любой цепи Маркова на произвольном фазовом пространстве, в т. ч. на произвольном топологическом пространстве, т. е. они применимы и к цепи Маркова, изучаемой в ~\cite{Cio2019}.\\

Конечно-аддитивные меры появились в операторной теории общих цепей Маркова около 60 лет тому назад. Однако, по сегодняшний день в литературе трудно найти конкретные (не вырожденные) и простые примеры использования таких мер в этой теории. 

Вместе с тем, конечно-аддитивные меры в теории цепей Маркова возникают не только при общих фазовых пространствах. Они могут дать кое-что новое даже для ``почти'' феллеровских ЦМ на отрезке $[0,1]$. 

В нашем материале, предварительно представленном в работе ~\cite[пункт 3, стр. 13--30]{Zhd06}
содержится подробный эргодический анализ пяти цепей Маркова, заданных на отрезке $[0,1]$, и имеющих всего лишь два возможных перехода из каждой точки $x\in(0,1)$. Доказывается существование инвариантных чисто конечно-аддитивных мер для этих ЦМ. При этом используются наши общие теоремы.

Не будем приводить эти примеры здесь, и отправляем интересующихся читателей к источнику ~\cite{Zhd06}.

\end{document}